\documentclass{article}
\usepackage{amsmath,amsfonts,amsthm,amssymb}
\usepackage{amscd}

    \advance \topmargin by -\headheight
    \advance \topmargin by -\headsep

    \evensidemargin \oddsidemargin
\newtheorem {theorem}    {Theorem}[section]

\newtheorem {lemma}      [theorem]    {Lemma}
\newtheorem {corollary}  [theorem]    {Corollary}
\newtheorem {proposition}[theorem]    {Proposition}
\newtheorem {claim}      [theorem]    {Claim}

\theoremstyle{definition}
\newtheorem*{definition}{Definition}
\theoremstyle{remark}
\newtheorem*{remark}{Remark}

\def\Bbb#1{\mathbb{#1}}

\def\ad{\operatorname{ad}}
\def\Aut{\operatorname{Aut}}
\def\depth{\operatorname{depth}}
\def\height{\operatorname{height}}

\def\char{{\rm char}}
\def\deg{{\rm deg}}
\def\LC{{\rm LC}}

\newcommand  {\QED}    {\def\qedsymbol{$\square$}\qed}
\def\Bcal{{\mathcal B}}
\def\Gcal{{\mathcal G}}
\def\Ncal{{\mathcal N}}

\def\Ucal{{\mathcal U}}

\def\Bhat{{\widehat{B}}}
\def\Ghat{{\widehat{G}}}

\def\Uhat{{\widehat{U}}}

\def\Bbar{{\overline{B}}}
\def\Gbar{{\overline{G}}}

\def\Ubar{{\overline{U}}}

\def\grg{{\mathfrak g}}
\def\grh{{\mathfrak h}}

\def\dbG{{\mathbb G}}

\def\dbN{{\mathbb N}}

\def\dbZ{{\mathbb Z}}

\def\eps{\varepsilon}

\def\lam{\lambda}            
\def\sig{\sigma}                
\def\phi{\varphi}

\def\Ker{{\rm Ker\,}}

\def\depth{{\rm depth}}

\def\wgal{\widetilde w}

\def\Gtil{\widetilde G}

\def\Wtil{\widetilde W}

\def\grg{{\mathfrak g}}
\def\la{\langle}
\def\ra{\rangle}
\def\T{\mathcal T}

\def\iff{\,\Leftrightarrow\,} 
\def\rank{\operatorname{rank}} 
\def\Z{{\mathbb Z}}
\def\negspace{\vspace{-0.2cm}}

\begin{document}

\author{Lisa Carbone\\
\footnotesize{Department of Mathematics, Hill Center, Busch Campus, Rutgers,
The State University of New Jersey}\\
{\footnotesize 110 Frelinghuysen Rd Piscataway, NJ 08854-8019} \\
\footnotesize{e-mail: carbonel@math.rutgers.edu} \\ \\ \\
Mikhail Ershov
\thanks{Corresponding author}
\\
\footnotesize{School of Mathematics, Institute for Advanced Study} \\
\footnotesize{1 Einstein Drive, Princeton, NJ 08540.} \\ \\
\footnotesize{Current address: University of Virginia, Department of Mathematics,}\\
\footnotesize{Kerchof Hall, Charlottesville, VA 22904}\\
\footnotesize{e-mail: ershov@virginia.edu}  \\ \\ \\
Gordon Ritter \\
\footnotesize{Department of Physics, Harvard University,} \\
\footnotesize{17 Oxford St., Cambridge, MA 02138} \\
\footnotesize{e-mail: ritter@post.harvard.edu}
}

\title{Abstract simplicity of complete Kac-Moody groups over finite fields}
\date{}
\maketitle

\noindent \thanks{\footnotesize The first author was supported in part by NSF grant \#DMS-0401107}

\noindent \thanks{\footnotesize The second author was supported by NSF grant \#DMS-0111298}
\footnote{This material is based in part upon work supported by the National Science
Foundation under agreement No. DMS-0111298. Any opinions, findings and
conclusions or recommendations expressed in this material are those of
the authors and do not necessarily reflect the views of the National Science
Foundation}

\noindent \thanks{2000 Mathematics subject classification. Primary 20E42; Secondary 20E32, 17B67, 20E18, 22F50}
\begin{abstract}
Let $G$ be a Kac-Moody group over a finite field corresponding to a generalized
Cartan matrix $A$, as constructed by Tits. It is known that $G$ admits the structure of a BN-pair, 
and acts on its corresponding building.
We study the complete Kac-Moody group $\Ghat$ which is defined to be the
closure of $G$ in the automorphism group
of its building. Our main goal is to determine when complete Kac-Moody groups
are abstractly simple, that is have no
proper non-trivial normal subgroups. Abstract simplicity of $\Ghat$ was
previously known to hold when A is of affine type. We extend this result to
many indefinite cases, including all hyperbolic generalized Cartan matrices $A$
of rank at least
four. Our proof uses Tits' simplicity theorem for groups with a BN-pair and
methods from the theory of pro-$p$ groups.
\end{abstract}

    \renewcommand{\baselinestretch}{1.2}

\section{Introduction}

Let $k$ denote a finite field which will remain fixed throughout the paper. Let $A$ be a generalized
Cartan matrix, and let $\dbG_{A}$ be the corresponding
Kac-Moody group functor of simply-connected type, as constructed
by Tits in \cite{Ti2}. The group $G(A):=\dbG_A(k)$  is usually called a
\emph{minimal} or \emph{incomplete} Kac-Moody group.

Distinct constructions of \emph{complete} Kac-Moody groups
are given in the papers of Carbone and Garland~\cite{CG} and
Remy and Ronan~\cite{RR}; in both cases the group constructed
is the completion of Tits' group $G(A)$ with respect to a certain
topology. In \cite{RR}, the topology comes from the action of $G(A)$
on its associated positive building; we denote the corresponding completion
by $\Ghat(A)$.\footnote{Such completions are called \it{topological Kac-Moody
groups }\rm in \cite{RR}.}
In \cite{CG}, one starts with the integrable highest-weight
module $V^{\lam}$ for the Kac-Moody algebra $\mathfrak g(A)$
corresponding
to a regular dominant integral weight $\lam$, then considers a certain
$\dbZ$-form, $V_{\dbZ}^{\lam}$ of $V^{\lam}$,
and the action of $G(A)$ on a $k$-form, $V_k^{\lam}$ of $V_{\dbZ}^{\lam}$.
The corresponding completion
of $G$ (which depends on $\lam$) will be denoted by $\widehat{G}^\lam(A)$.
We describe these constructions in detail in Section 2.

Though constructed in different ways, the groups $\Ghat(A)$ and $\Ghat^{\lam}(A)$
have very similar structure.
In particular, in both cases the complete Kac-Moody group 
is locally compact and totally disconnected. Recently, Baumgartner and Remy
showed\footnote{Private communication.}
that for any weight $\lam$, the Remy-Ronan completion $\Ghat(A)$
is a homomorphic image of the Carbone-Garland completion $\widehat{G}^\lam(A)$
(see Theorem~\ref{relation} for a
precise statement).

The main goal of this paper is to investigate
when complete Kac-Moody groups are abstractly simple.
Recall that a topological group is called {\it abstractly simple }
if it has no proper non-trivial normal subgroups, and
{\it topologically simple } if it has no proper non-trivial closed normal subgroups.
In view of Theorem~\ref{relation}, the simplicity question should be asked for the smaller groups $\Ghat(A)$.
An obvious necessary condition for the simplicity of $\Ghat(A)$ is that $A$ should be indecomposable:
if $A_1, \ldots, A_k$ are indecomposable blocks of $A$, then $\Ghat(A)\cong \prod_{i=1}^k \Ghat(A_i)$.
Remy~\cite[Theorem 2.A.1]{Re2} proved\footnote{Remy considers a more general class of quasi-split Kac-Moody groups.}
that if $A$ is indecomposable, then
the group $\Ghat(A)$ is topologically simple when $|k|>3$,
and asked whether abstract simplicity holds as well \cite[Question 30]{Re3}.
We answer this question in the affirmative for a large
class of Kac-Moody groups.

\begin{theorem}
Let $A$ be an indecomposable generalized Cartan matrix.
Assume that one of the following holds:\negspace
\begin{enumerate}
\item[(a)] $|k|>3$, $p=\char(k)>2$, $A$ is symmetric, and any $2\times 2$ submatrix of $A$ is of finite or 
affine type;
\item[(b)] $|k|>3$ and any $2\times 2$ submatrix of $A$ is of finite type.
\end{enumerate}
\negspace
Then the group $\Ghat(A)$ is abstractly simple.
\label{main}
\end{theorem}

\begin{remark}
In the hypotheses of Theorem~\ref{main}, a submatrix is not necessarily proper.
\end{remark}

Abstract simplicity of $\Ghat(A)$ was previously known only when $A$ is of finite type
(in which case $\Ghat(A)$ is a finite group) or $A$ is of affine type, in which case $\Ghat(A)$
is isomorphic to the group of $K$-points of a simple algebraic group defined over $K=k((t))$
\cite{Re2}.
The class of matrices covered by Theorem~\ref{main} includes many indefinite examples,
including all hyperbolic matrices of rank at least four -- see Proposition~\ref{affine_hyperbolic}.

\begin{remark}
Recently, Caprace and Remy~\cite{CR} proved (abstract)
simplicity of the incomplete group $G(A)$ modulo its finite center in the case
when the associated Coxeter group is not affine and assuming that $|k|$ is
sufficiently large.
If $A$ is affine, the incomplete group $G(A)$ (modulo its center) is infinite and residually finite,
and hence cannot be simple.
\end{remark}

Briefly, our approach to proving  Theorem~\ref{main} will be as follows.
A celebrated theorem of Tits \cite{B} gives sufficient conditions (called ``simplicity axioms")
for a group with a BN-pair to be simple. The group $\Ghat(A)$
routinely satisfies most of these axioms if $|k|>3$ (for arbitrary $A$),
which already implies topological simplicity of $\Ghat(A)$. It is not clear
if the remaining axioms hold in general; however, they do hold if the ``positive
unipotent'' subgroup $\Uhat(A)$ is topologically finitely generated -- see
Theorem~\ref{simprel}
($\Uhat(A)$ is defined as the completion of the group $U(A)$ generated by all
positive root subgroups).
This follows from the fact that $\Uhat(A)$ is a pro-$p$ group and basic
properties of pro-$p$ groups (see \cite[Chapter 1]{DSMS}).

Thus our main task is to prove (topological) finite generation of $\Uhat(A)$
under the hypotheses of Theorem~\ref{main}.
If any $2\times 2$ submatrix of $A$ is of finite type and $|k|>3$,
then $U(A)$ is finitely generated (as an abstract group) by a theorem of
Abramenko~\cite{A}, which immediately implies finite generation of
$\Uhat(A)$. This proves Theorem~\ref{main}(b).

In Section~6, we show (see Theorem~\ref{fingen}) that $\Uhat(A)$
is finitely generated as long as 

(a) $A$ is symmetric;

(b) for every $2\times 2$ submatrix $C$ of $A$, the group $\Uhat(C)$ is ``well behaved''.

The latter is a certain technical condition which is ``almost''
equivalent to generation by simple root subgroups.
In Section~7, we show that the group $\Uhat(C)$ is well behaved
when $C$ is a $2\times 2$ matrix of finite or affine type (assuming $p>3$), 
using an explicit realization of $\Uhat(C)$. This completes the proof of Theorem~\ref{main}(a).
When $C$ is a $2\times 2$ hyperbolic matrix, we do not know if $\Uhat(C)$
is well behaved or even if  $\Uhat(C)$ is finitely generated. It seems
that essentially new ideas are needed to settle this case.

\bf{Acknowledgements. }\rm
We are extremely grateful to Bertrand Remy and Pierre-Emmanuel Caprace
for  illuminating discussions which led to significant
improvement of the results in this paper. We are indebted to an anonymous referee
for pointing out a mistake in the earlier version of this paper.
We wish to thank Bertrand Remy
for his careful reading of this paper, and for many useful suggestions.
We thank Howard Garland who suggested this project and who                                                                                                   
contributed to the early stages of it.
We are also grateful to Udo Baumgartner for helpful discussions,
and to Jacqui Ramagge for providing us with a copy of her thesis~\cite{Ra1}.
\vskip .15cm
\bf{Some conventions about topological groups. }\rm
Let $G$ be a topological group. By a base (resp. subbase) for the topology of $G$
we will mean a base (resp. subbase) of neighborhoods of the identity.
Recall that a topological group $G$ is {\it topologically generated by a set $S$}
if the subgroup abstractly generated by $S$ is dense in $G$.
When discussing finite generation of pro-$p$ groups, we will always
be interested in topological generating sets, so the word `topological' will often be omitted.
Since (infinite) pro-$p$ groups are never finitely generated as abstract groups,
this convention should not cause any confusion.

\section{Kac-Moody algebras and groups}
\label{sec:prelim}

\subsection{Generalized Cartan matrices}

Let $I=\{1,2,\dots,l\}$ be a finite set.
A matrix $A=(a_{ij})_{i,j\in I}$ is called a \it{generalized
Cartan matrix }\rm if its entries satisfy the following conditions:
\begin{align*}
&\mbox{(a)}\quad a_{ij}\in{\Bbb Z} \mbox{ for } i,j\in I, &
&\mbox{(b)}\quad a_{ii}=2, \mbox{ for } i\in I&\\
&\mbox{(c)}\quad a_{ij}\leq 0 \mbox{  if } i\neq j&
&\mbox{(d)}\quad a_{ij}=0 \iff a_{ji}=0.&\\
\end{align*}
\vskip -.4cm
By a \emph{submatrix} of $A$, we mean a matrix of the form
$$A_{J}=(a_{ij})_{i,j\in J},$$ where $J$ is a subset of $I$.
We say that the submatrix $A_J$ is \it{proper }\rm if $J\neq I$.
The matrix $A$ is called {\it indecomposable} if there is no partition of
the set $I$ into two non-empty subsets so that
$a_{ij}=0$ whenever $i$ belongs to the first subset, while $j$
belongs to the second. A submatrix $A_J$ is called an \it{indecomposable block }\rm
of $A$ if $A_J$ is indecomposable, and $J$ is maximal with this property.

A generalized Cartan matrix $A$ is said to be of
\begin{enumerate}
\item[]
{\it finite (classical) type} if $A$ is positive definite;
in this case, $A$ is the Cartan matrix of a finite dimensional
semisimple Lie algebra,
\item[]
{\it affine type} if $A$ is positive semi-definite, but not
positive definite,
\item[]
{\it indefinite type} if $A$ is neither of finite nor affine type.
\end{enumerate}

An indefinite matrix is said to be of \it{hyperbolic type }\rm 
(in the sense of \cite[5.10, p.66]{Kac}) if every proper
indecomposable submatrix of $A$ is of finite or affine type. 

\begin{proposition}
\label{affine_hyperbolic}
Let $A$ be an indecomposable $l\times l$ matrix, and
let $B$ be an $s\times s$ submatrix of $A$. Assume that either 
\begin{itemize}
\item[(i)] $A$ is of finite type,
\item[(ii)] $A$ is of affine type and $s\leq l-1$ (that is, $B$ is a proper submatrix of $A$),
\item[(iii)] $A$ is of hyperbolic type and $s\leq l-2$.
\end{itemize}
Then $B$ is of finite type.
\end{proposition}
\begin{proof}
Case (i) is obvious, and case (ii) is well known  \cite[Chapter~4]{Kac}. 
Now let $A$ be hyperbolic and $s\leq l-2$. 
We can assume that $B$ is indecomposable, since if
every indecomposable block of $B$ is of finite type, then so is $B$.
Since $A$ is indecomposable, there exists an indecomposable $(s+1)\times (s+1)$
submatrix $C$ of $A$ such that $C$ contains $B$. Since $s+1<l$,
the matrix $C$ must be of finite or affine type, whence $B$ is of finite type
by cases (i) and (ii).
\end{proof}

\subsection{Kac-Moody algebras}

For the rest of this section, fix a generalized Cartan matrix $A=(a_{ij})_{i,j\in I}$,
and let $l=|I|$. A \it{realization }\rm of $A$ over $\mathbb Q$ is a triple
$(\grh, \Pi, \Pi^{\vee})$ where $\grh$ is a vector space over $\mathbb Q$
of dimension \mbox{$2l-\rank(A)$}, and $\Pi=\{\alpha_1,\dots,\alpha_l\}\subseteq
\mathfrak{h}^*$
and $\Pi^{\vee}=\{\alpha_1^{\vee},\dots,\alpha_l^{\vee}\}\subseteq
\mathfrak{h}$ are linearly independent sets, such that
$\la\alpha_j,\alpha_i^{\vee}\ra=a_{ij}$ for $i,j\in I$. As usual, $\la\cdot ,\cdot \ra$
denotes the natural pairing between $\mathfrak{h}$ and $\mathfrak{h}^*$.
Elements of $\Pi$ are called {\it simple roots}
and elements of $\Pi^{\vee}$ {\it simple coroots}.

The associated  Kac-Moody algebra $\grg=\grg_{A}$ is a Lie algebra over
${\mathbb Q}$, generated by $\mathfrak{h}$ and elements $(e_i)_{i\in I}$,
$(f_i)_{i\in I}$ subject to the Serre-Kac relations:
\begin{align*}
&\mbox{1) }\ [\mathfrak{h}, \mathfrak{h}]=0&
&\mbox{5) }\ [e_i,f_j]=0,\ i\neq j&\\
&\mbox{2) }\ [h,e_i]=\la\alpha_i,h\ra e_i, h\in \mathfrak{h}&
&\mbox{6) }\ (\ad\, e_i)^{-a_{ij}+1}(e_j)=0,\ i\neq j&\\
&\mbox{3) }\ [h,f_i]=-\la\alpha_i, h\ra f_i, h\in \mathfrak{h}&
&\mbox{7) }\ (\ad\, f_i)^{-a_{ij}+1}(f_j)=0,\ i\neq j&\\
&\mbox{4) }\ [e_i,f_i]=\alpha_i^{\vee}&
&&\\
\end{align*}
It is easy to see that $\grg_A$ depends only on $A$ and not on its realization
(see \cite[1.1]{Kac}).

Relative to $\mathfrak{h}$, the Lie algebra $\mathfrak{g}$ has decomposition
$\mathfrak{g}=\mathfrak{h}\oplus\bigsqcup_{\alpha\in\Delta} \mathfrak{g}^{\alpha},$ where
$$\mathfrak{g}^{\alpha}\quad=\quad\{x\in \mathfrak{g}\mid [h,x]=\la\alpha,h\ra x,\ h\in\mathfrak{h}\},$$
and $\Delta=\{\alpha\in \mathfrak{h}^* \backslash\{0\} \mid
\mathfrak{g}^{\alpha} \neq 0 \}$. Elements of $\Delta$ are called the \it{roots }\rm of $\mathfrak{g}$.
Each root has the form $\sum_{i\in I} n_i\alpha_i$ where $n_i\in\Z$ and either
$n_i\geq 0$ for all $i$, or $n_i\leq 0$ for all $i$. The roots are called positive or negative
accordingly; the set of positive (resp.~negative) roots will be denoted by
$\Delta^+$ (resp.~$\Delta^-$).
The \it{height }\rm of a root $\alpha=\sum_{i=1}^l n_i\alpha_i$ is defined
to be the integer $\sum_i n_i$.

\subsection{Real roots and the Weyl group}

For $i\in I$ define $w_i\in \Aut(\mathfrak h^*)$ by setting
$w_i(\alpha)= \alpha-\la\alpha,\alpha_i^{\vee}\ra\alpha_i$.
The group $W = \langle \{w_i\} \rangle$ generated by the $w_i$ is called the
\emph{Weyl group} associated to $A$.
The set $\Phi= W(\Pi)$ is a subset of $\Delta$, called the set of \emph{real
roots.} The remaining roots $\Delta\backslash \Phi$ are called
\emph{imaginary roots.}

The Weyl group $W$  has a faithful action on $\grh$ defined by
$w_i(h)= h-\la\alpha_i,h\ra\alpha_i^{\vee}$. Moreover, the pairing
$\la\cdot,\cdot\ra$ is
$W$-invariant, that is, $\la w\alpha, wh \ra=\la \alpha,h\ra$ for
$\alpha\in\grh^*$, $h\in\grh$ and $w\in W$.

For each real root $\alpha$, define the corresponding coroot $\alpha^{\vee}$ as follows:
write $\alpha$ in the form $w\alpha_i$ for some $w\in W$ and $i\in I$
and set $\alpha^{\vee}=w\alpha_i^{\vee}$. One can show (see \cite{Kac}) that
$\alpha^{\vee}$ is independent of the above choice.
The correspondence $\alpha\mapsto\alpha^{\vee}$ is not linear; however,
it does satisfy some nice properties:
\begin{proposition}
The following hold:\negspace
\begin{enumerate}
\item[(a)] For each $\alpha\in\Phi$, the coroot $\alpha^{\vee}$ is an integral
linear combination of
$\{\alpha_i^{\vee}\}$, and the coefficients are all non-negative (resp. non-positive) if $\alpha\in\Phi^+$
(resp. $\alpha\in\Phi^-$). Furthermore,
\mbox{$(-\alpha)^{\vee}=-\alpha^{\vee}$}.

\item[(b)] Given $\alpha,\beta\in\Phi$, we have
$\la\alpha,\beta^{\vee}\ra >0$ (resp. $=0, <0$) $\iff$ $\la\beta,\alpha^{\vee}\ra >0$ (resp. $=0, <0$).

\item[(c)] For every $\alpha\in\Phi$, we have $\la\alpha,\alpha^{\vee}\ra=2$.
\end{enumerate}
\end{proposition}

In order to define Kac-Moody groups, we introduce a related group $W^*\subseteq \Aut (\grg)$.
By definition, $W^*$ is generated by elements $\{w_i^*\}_{i \in I}$, where
\[
	w_i^*
	=
	\exp(\ad e_i)\exp(-\ad f_i)\exp(\ad e_i)
	=
	\exp(-\ad f_i)\exp(\ad e_i)\exp(-\ad f_i).
\]
The group $W^*$ is a central extension of $W$.
More specifically, there is a surjective homomorphism $\eps: W^*\to W$
which sends $w_i^*$ to $w_i$ for all $i$; the kernel of $\eps$ is an elementary
abelian group of exponent $2$ generated by $\{(w_i^*)^2\}$, as  follows
immediately from \cite[3.3]{Ti2}.

Finally, we define certain elements $\{e_\alpha\in\mathfrak g\}_{\alpha\in \Phi}$.
Given $\alpha\in \Phi$, write $\alpha$ in the form $w\alpha_j$ for some $j\in I$ and $w\in W$,
choose $w^*\in W^*$ which maps onto $w$, and set $e_{\alpha}=w^* e_{\alpha_j}$.
It is clear from \cite[(3.3.2)]{Ti2} that $e_{\alpha}$ belongs to $g^{\alpha}$,
$e_{\alpha}$ is uniquely determined up to  sign, and for all $i \in I$,
$w_i^* e_{\alpha}=\eta_{\alpha,i}e_{w_i\alpha}$ for some
constants $\eta_{\alpha,i} \in \{ \pm 1 \}$.
These constants $\{\eta_{\alpha,i}\}$ will appear in the definition of Kac-Moody
groups.

\subsection{Kac-Moody groups and Tits' presentation}

The construction of (incomplete) Kac-Moody groups over arbitrary fields
is due to Tits \cite{Ti2}. One may define these
groups by generators and relators. While not explicitly stated in Tits' paper,
such a presentation appears in the papers of Carter~\cite{C} and (in a slightly
different form) Morita and Rehmann \cite{MR}.

The group $G=G(A)$ defined below is called the \emph{incomplete simply-connected
Kac-Moody group} corresponding to $A$.
The presentation we use is  ``almost canonical'' except for the choice
of elements $\{e_{\alpha}\}$ which determine the constants $\{\eta_{\alpha,i}\}$.

By definition, $G(A)$ is generated by the set of symbols
$\{\chi_{\alpha}(u)\mid \alpha\in \Phi, u\in k\}$ satisfying
relations (R1)-(R7) below. In all the relations
$i,j$ are elements of $I$,
$\,u,v$ are elements of $k$ (arbitrary, unless mentioned otherwise)
and $\alpha$ and $\beta$ are real roots.

(R1) $\chi_{\alpha}(u+v)=\chi_{\alpha}(u)\chi_{\alpha}(v)$;

(R2) Let $(\alpha,\beta)$ be a \it{prenilpotent pair, }\rm
that is, there exist $w,\ w'\in W$ such that
$$w\alpha, \ w\beta\in\Phi^+{\text{ and }}w'\alpha, \ w'\beta\in\Phi^-.$$
Then
\[
	[\chi_{\alpha}(u),\chi_{\beta}(v)]
	=
	\prod_{m,n \geq 1}
	\chi_{m\alpha+n\beta}(C_{mn\alpha\beta}u^m v^n)
\]
where the product on the right-hand side is taken over all real roots
of the form $m\alpha+n\beta$, $m,n\geq 1$, in some fixed order, and
$C_{mn\alpha\beta}$ are integers independent of $k$ (but depending on the
order).

\noindent For each $i\in I$ and $u\in k^*$ set

$\chi_{\pm i}(u)=\chi_{\pm \alpha_i}(u)$,

$\wgal_{i}(u)=\chi_{i}(u)\chi_{-i}(-u^{-1})\chi_{i}(u)$,

$\wgal_{i}=\wgal_{i}(1)$ and $h_{i}(u)=\wgal_{i}(u)\wgal_{i}^{-1}$.

\noindent The remaining relations are

(R3) $\wgal_{i}\chi_{\alpha}(u)\wgal_{i}^{-1}=
\chi_{w_{i}\alpha}(\eta_{\alpha,i}u)$,

(R4)
$h_{i}(u)\chi_{\alpha}(v)h_{i}(u)^{-1}=\chi_{\alpha}(vu^{\la\alpha,\alpha_i^{
\vee}\ra})$ for $u\in k^*$,

(R5) $\wgal_{i}h_j(u)\wgal_i^{-1}=h_j(u)h_i(u^{-a_{ji}})$,

(R6) $h_{i}(uv)=h_{i}(u)h_{i}(v)$ for $u,v\in k^*$, and

(R7) $[h_{i}(u),h_{j}(v)]=1$ for $u,v\in k^*$.

\noindent An immediate consequence of relations (R3) is that $G(A)$
is generated by $\{\chi_{\pm i}(u)\}$.

\begin{remark}
The product on the right-hand side of (R2) is finite because
of the following well-known criterion:
a pair $\{\alpha,\beta\}$ is prenilpotent
if and only if $\alpha\neq -\beta$ and $|(\Z_{>0}\alpha +\mathbb
Z_{>0}\beta)\cap \Phi|<\infty$. This result is proved in \cite[Proposition 4.7]{KP}
for $\alpha,\beta\in \Phi^+$, and the general case follows from
the straightforward observation that for any $\alpha,\beta\in \Phi$
with $\alpha\neq -\beta$ there exists $w\in W$ such that either
$w\alpha,\,w\beta\in\Phi^+$ or $w\alpha,\,w\beta\in\Phi^-$.
\end{remark} 

\vskip .1cm
Intuitively, one should think of the above presentation as an analogue of the
Steinberg presentation for
classical groups with $\chi_{\alpha}(u)$ playing the role of
$\exp(ue_{\alpha})$. In the next
subsection we give a representation-theoretic interpretation of Kac-Moody
groups which makes the above analogy precise.

Next we introduce several subgroups of $G = G(A)$:
\begin{enumerate}
\item {\it Root subgroups $U_{\alpha}$}. 
For each $\alpha\in\Phi$ let $U_{\alpha}=\{\chi_{\alpha}(u)\mid u\in k\}$.
By relations (R1), each $U_{\alpha}$ is isomorphic to the additive
group of $k$.

\item {\it The ``extended" Weyl group $\Wtil$}. Let $\Wtil$ be the subgroup of
$G$ generated by elements
$\{\wgal_{i}\}_{i\in I}$. One can show that $\Wtil$ is isomorphic to the group $W^*$
introduced before, so there is a surjective homomorphism
$\eps:\widetilde W\to W$ such that $\eps(\wgal_{i})=w_{i}$ for $i\in I$. Given
$\wgal\in \Wtil$ and $w\in W$, we will say that $\wgal$ is a representative of
$w$ if $\eps(\wgal)=w$.  It will be convenient to identify (non-canonically) $W$
with a subset (not a subgroup) of $\Wtil$
which contains exactly one representative of every element of $W$. By abuse of
notation, the set of those representatives will also be denoted by $W$.
It follows from relations (R3) that $w U_{\alpha} w^{-1}=U_{w\alpha}$ for any
$\alpha\in \Phi$ and $w\in W$.

\item {\it ``Unipotent'' subgroups.} Let  $U^{+}=\la U_{\alpha}\mid \alpha\in
\Phi^+\ra$, and $U^{-}=\la U_{\alpha}\mid \alpha\in \Phi^-\ra$.

\item {\it ``Torus'' (``diagonal'' subgroup).} Let $H=\la\{h_{i}(u)\mid i\in I, u\in k\}\ra$.
One can show that relations (R6)-(R7) are defining relations for $H$, so
$H$ is isomorphic to the direct sum of $l$ copies of $k^*$.

\item {\it ``Borel'' subgroups.} Let $B^+=\la U^+, H\ra$ and $B^-=\la U^-, H\ra$.
By relations (R4), $H$ normalizes both $U^{+}$ and $U^{-}$, so we have
$B^+=H U^{+}=U^+ H$ and $B^-=H U^{-}=U^- H$.

\item {\it ``Normalizer.''} Let $N$ be the subgroup generated by
$\widetilde{W}$ and $H$.
Since $\widetilde{W}$ normalizes $H$, we have $N=\widetilde W H$.
It is also easy to see that $N/H\cong W$.
\end{enumerate}

Tits~\cite{Ti2} proved that $(B^+,N)$ and $(B^-,N)$ are
BN-pairs\footnote{Recall that BN-pairs are also called Tits systems.} of $G$.
In fact, $G$ admits the stronger structure of a twin BN-pair, but we will not use it. 
Let $X^+$ and $X^-$ be the buildings associated with $(B^+,N)$ and $(B^-,N)$, respectively.
Since the field $k$ is finite, the buildings
$X^+$ and $X^-$ are locally finite as chamber complexes. In fact, $X^+$
and $X^-$ have constant thickness $|k|+1$ (see \cite[Appendix KMT]{DJ}).

Below we list some of the fundamental properties of these BN-pairs.
We will work mostly with the positive BN-pair $(B^+,N)$, and from now on, write
$B$ for $B^+$ and $U$ for $U^+$.

\negspace
\begin{enumerate}
\item[(a)] $B\cap N=H$, so the Coxeter group associated to $(B,N)$ is isomorphic to the Weyl group $W=W(A)$;

\item[(b)] Bruhat decomposition: $G=BWB$;

\item[(c)] Birkhoff decomposition: $G=U^- W B= B^- W U= U W B^- = BW U^-$.
\end{enumerate}
\negspace

Of course, (b) follows directly from $(B,N)$ being a BN-pair, and the proof of (c) uses the twin
BN-pair structure (see \cite{KP}).

Finally, we shall need a presentation by generators and relators for the
group $U$ established by Tits~\cite[Proposition 5]{Ti1}.
\begin{theorem}
\label{pres_U}
The group $U$ is generated by the elements $\{\chi_{\alpha}(u)\mid \alpha\in\Phi^+, u\in k\}$
subject to relations (R1) and (R2) defined earlier in this section.
\end{theorem}

\subsection{Representation-theoretic interpretation of Kac-Moody groups}

The following interpretation of Kac-Moody groups was given by Carbone
and Garland \cite{CG} (see also \cite{Ti3}). This construction generalizes
that of Chevalley groups \cite{St}.
Let ${\mathcal U}$ be the universal enveloping algebra of
$\mathfrak{g}$. Let $\Lambda\subseteq \mathfrak{h}^*$  be the linear span of  $\alpha_i$, for $i\in I$, and
$\Lambda^{\vee}\subseteq
\mathfrak{h}$  be the linear span of
$\alpha_i^{\vee}$, for $i\in I$.
Let
${\mathcal U}_\dbZ\subseteq {\mathcal U}$ be the $\dbZ$-subalgebra generated by
$e_i^m/m!$, $f_i^m/m!$, and $\binom{h}{m}$, for $i\in I$,  $h\in\Lambda^{\vee}$ and $m\geq 0$. Then
${\mathcal U}_\dbZ$ is a $\dbZ$-form of ${\mathcal U}$, i.e.  ${\mathcal U}_\dbZ$ is a subring and the canonical map
${\mathcal U}_{\dbZ}\otimes{\mathbb Q}\longrightarrow {\mathcal U}$ is bijective.
For a field $K$, let ${\mathcal U}_{K}={\mathcal U}_\dbZ \otimes K$, and $\mathfrak{g}_{K}=\mathfrak{g}_\dbZ \otimes K$.

Now let $\lambda\in \mathfrak{h}^{*}$ be a regular dominant integral weight,
that is, $\la\lambda, \alpha_i^{\vee}\ra\in \dbZ_{>0}$ for every $i\in I$.
Let $V^{\lambda}$ be the
corresponding irreducible highest weight module. Choose
a highest-weight vector $v_{\lam}\in V^{\lam}$, and let
$V^{\lambda}_\dbZ \subset V^{\lam}$
be the orbit of $v_{\lam}$ under the action of ${\mathcal U}_\dbZ$.
Then $V^{\lambda}_\dbZ$ is a $\dbZ$-form of $V_{\lam}$ as well as a
${\mathcal U}_\dbZ$-module.
Similarly, $V^{\lambda}_k:= k\otimes_\dbZ V^{\lambda}_\dbZ$
is a ${\mathcal U}_{k}$-module.

It is straightforward to establish the following (see \cite[Proposition 3]{Ti3}).
\begin{proposition}
There is a (unique) homomorphism $\pi_{\lam}:G\to \Aut(V^{\lambda}_k)$ such that
\begin{align*}
&\pi_{\lam}(\chi_{\alpha_i}(u))\quad=\quad\sum_{m=0}^{\infty} u^m \dfrac{e_i^m}{m!}
\quad\mbox{ for } i\in I \mbox{ and } u\in k,&\\
&\pi_{\lam}(\chi_{-\alpha_i}(u))\quad=\quad\sum_{m=0}^{\infty} u^m \dfrac{f_i^m}{m!}
\quad\mbox{ for } i\in I \mbox{ and } u\in k.&
\end{align*}
\label{incompletecomparison}
\end{proposition}
\vskip -.4cm
The expressions on the right-hand side are well defined
automorphisms of $V^{\lambda}_k$ since $e_i$ and $f_i$ are locally nilpotent on
$V^{\lambda}_k$. Let $G^{\lam}=\pi_{\lam}(G)$. As we will see later in this section,
the kernel of $\pi_{\lam}$ is finite, central and contained in $H$.

\subsection{Complete Kac-Moody groups}

As mentioned in the introduction, distinct completions of $G(A)$ were
given in the papers of Carbone and Garland~\cite{CG} and Remy and
Ronan~\cite{RR}. We now briefly review
these constructions, starting with the Remy-Ronan completion.
As above, let $X^+$ be the building
associated with the positive BN-pair $(B,N)$, and
consider the action of $G$ on $X^+$.
Recall that $X^+$ is locally finite as a chamber complex.
Define the topology on $G$
by the subbase\footnote{Recall that by `subbase' for a topology
on a group, we mean a `subbase of neighborhoods of the identity'.}
consisting of stabilizers of vertices of $X^+$ or, equivalently,
fixators (pointwise stabilizers) of chambers of $X^+$.
We shall call this topology the \emph{building topology}.
The completion of $G$ in its building topology will be referred
to as the \emph{Remy-Ronan completion} and denoted by $\Ghat$.

We will  make few references to  the action of $G$
on its building in this work. All we will need is the description of the building
topology in purely group-theoretic terms. Since $(B,N)$ is a BN-pair,
we know that \negspace
\begin{enumerate}
\item[(a)] The subgroup $B$ of $G$ is a chamber fixator,

\item[(b)] $G$ acts transitively on the set of chambers of $X^+$.
\end{enumerate}
\negspace Therefore, the family $\{gB g^{-1}\}_{g\in G}$
is a subbase for the building topology.

Let $Z$ be the kernel of the natural map $G\to \Ghat$
(or, equivalently, the kernel of the action of $G$ on $X^+$).
Using results of Kac and Peterson \cite{KP}, Remy and Ronan~\cite[1.B]{RR} showed
that $Z$ is a subgroup of $H$ (and hence finite); furthermore,
$Z$ coincides with the center of $G$.

Now let $\Bhat$ (resp. $\Uhat$) be the closure
of $B$ (resp. $U$) in $\Ghat$. The natural images of $N$ and $H$
in $\Ghat$ are discrete, and therefore we will denote
them by the same symbols (without hats). This involves
some abuse of notation since the image of $H$ in
$\Ghat$ is isomorphic to $H/Z$.

The following theorem is a collection of results from \cite{Re1} and \cite{RR}:

\begin{theorem} \label{collection}
Let $\Ghat$, $\Bhat$ and $N$ be as above. The following hold:
\negspace
\begin{enumerate}
\item[(a)] The pair $(\Bhat,N)$ is a BN-pair of $\Ghat$. Moreover,
if $\widehat X^+$ is the associated building, there exists a $\Ghat$-equivariant
isomorphism between $X^+$ and $\widehat X^+$.
In particular, the Coxeter group associated to $(\Bhat,N)$
is isomorphic to $W=W(A)$.

\item[(b)] The group $\Bhat$ is an open profinite subgroup of $\Ghat$. Furthermore,
$\Uhat$ is an open pro-$p$ subgroup of $\Bhat$.
\end{enumerate}
\end{theorem}

Now we turn to the Carbone-Garland completion.
Let $\lam$ be a regular weight, and let
$G^{\lam}=G^{\lam}(A)$ and $V^{\lam}_{k}$ be defined as in the previous subsection.
Now we define the \emph{weight topology} on $G^{\lam}$
by taking stabilizers of elements of $V^{\lam}_{k}$ as a subbase of
neighborhoods of the identity. The completion of $G^{\lam}$
in this topology will be referred to as the \emph{Carbone-Garland completion}
and denoted by $\Ghat^{\lam}(A)$. Since $G^{\lam}(A)$ is a homomorphic image of
$G(A)$, we can think of $\Ghat^{\lam}(A)$ as a completion of $G$ (and not $G^{\lam}$).
Let $\Bhat^{\lam}$ (resp. $\Uhat^{\lam}$) be the closures of $B$ (resp. $U$)
in $\Ghat^{\lam}(A)$. Then the obvious analogue of Theorem~\ref{collection}
holds; the fact that $\Uhat^{\lam}$ is a pro-$p$ group will be proved at the end
of this section (see Proposition~\ref{prop_CarboneGarland}); for all other
assertions see \cite[Section~6]{CG}.

The following relationship between the Remy-Ronan and Carbone-Garland completions
was established\footnote{Private communication.} by Baumgartner and Remy:

\begin{theorem}
For any regular weight $\lam$, there exists a (canonical)
continuous surjective homomorphism $\eps_{\lam}:\widehat{G}^\lam\to\Ghat$.
The kernel $K_{\lam}$ of $\eps_{\lam}$ is equal to $\bigcap\limits_{g\in\Ghat}g\Bhat^{\lam}g^{-1}$.
\label{relation}
\end{theorem}

It follows from Theorem~\ref{relation} that the kernel of the map $\pi_{\lam}: G\to G^{\lam}$
is finite and central. Indeed, consider the sequence of homomorphisms
\[
	\begin{CD}
	G @>\pi_\lam >> G^\lam @> >> \Ghat^\lam @> \eps_\lam >> \Ghat
	\end{CD}
\]
Clearly, the composition of these three maps is the natural map from $G$ to $\Ghat$.
We know that the kernel of the latter map is finite and central, hence the same should be true
for $\pi_{\lam}$.

In the case when $A$ is an affine matrix, Garland \cite{G} showed that
$K_{\lam}$ is a central subgroup of $H\subseteq \Bhat^{\lam}$ (and hence
finite). It is not clear to us how large $K_{\lam}$ can be in general.
Since $\Bhat^{\lam}$ is a profinite group, so is $K_{\lam}$;
furthermore, $K_{\lam}$ has a finite index pro-$p$ subgroup, which
follows from Proposition~\ref{prop_CarboneGarland} below.

\begin{proposition}
Let $\Uhat^{\lam}$ be the closure of $U$ in $\Ghat^{\lam}$.
The group $\Uhat^{\lam}$ is a pro-$p$ group.
\label{prop_CarboneGarland}
\end{proposition}
\proof
In \cite{CG}, it is shown that the $k$-vector space
$V^{\lambda}_{k}$ admits a basis $\Psi=\{v_1,v_2,v_3,\ldots \}$
consisting of weight vectors, that is, for each $i\in\dbN$
there exists a weight $\mu_i$ of $V^{\lambda}$ such
that $v_i$ lies in the weight component $V^{\lambda}_{\mu_i}$.
Each weight $\mu$ of $V^{\lambda}$ is of the form $\mu=\lambda-\sum_{i=1}^l k_i\alpha_i$,
where $k_i\in{\Bbb Z}_{\geq 0}$. Define the {\it depth} of
$\mu$ to be $\depth(\mu)=\sum_{i=1}^l k_i.$
For convenience, we order the elements of $\Psi$ such that
$\depth(\mu_i)\leq \depth(\mu_j)$ if $i<j$.

For each $n\geq 1$ let $V_n$ be the $k$-span of the set $\{v_1, v_2,\ldots, v_n\}$.
The group $U$ stabilizes $V_n$; moreover, it acts by upper-unitriangular
matrices (with respect to the above basis). Therefore,
we have a homomorphism $\pi_n:U\to GL_n(k)$ whose image
is a finite $p$-group (since $k$ has characteristic $p$).
Then $U_n = \Ker \pi_n$ consists of elements of $U$ which fix $V_n$ pointwise.
Since \mbox{$\bigcup_{n\geq 1}V_n = V^{\lambda}_{k}$}, the groups
$\{U_n\}_{n=1}^{\infty}$
form a base for the weight topology on $U$. Since each $U_n$ is a normal subgroup
of $U$ of $p$-power index, the completion of $U$ with respect to the weight topology
is a pro-$p$ group.
\QED

\begin{remark}
Remy and Ronan~\cite{RR} prove that
$\Uhat$ is a pro-$p$ group using its action on the building of $\Ghat$.
This fact can also be deduced from
Propostion~\ref{prop_CarboneGarland} by applying Theorem~\ref{relation}.
\end{remark}

We finish this section by describing explicitly the groups $G(A)$, $\Ghat(A)$
and $\Ghat^{\lam}(A)$
in the special case $A=A_{d-1}^{(1)}$ for some $d\geq 2$ (in the notation of
\cite[Chapter~4]{Kac});
an analogous result holds for any affine matrix $A$ -- see \cite{G}.
In this case, the incomplete group $G(A)$ is isomorphic to a central extension
of the group $SL_d(k[t,t^{-1}])$
by $k^*$. The Remy-Ronan completion $\Ghat(A)$ is isomorphic to $PSL_d(k((t)))$.
It is easy to see that the building topology on $\Ghat(A)$ coincides with the topology
on $PSL_d(k((t)))$ induced from the local field $k((t))$.

Recall (see above) that the center of $G(A)$ always lies in
the kernel of the natural map $G(A)\to\Ghat(A)$. On the other hand, the center of $G(A)$
does usually have non-trivial image in the Carbone-Garland completion $\Ghat^{\lam}(A)$.
For any $\lam$, there is a commutative diagram
\[
\begin{CD}
1@> >>k^*@> >>\widehat{SL_d}\Bigl(k((t))\Bigr)@> >>SL_d\Bigl(k((t))\Bigr)@> >>1\\
@.@.@VV \rho_{\lam}V @VV \rho V\\
{}@.{}@.\Ghat^{\lam}(A)@>\eps_{\lam} >> \Ghat(A)
\end{CD}
\]
where $\widehat{SL_d}(k((t)))$ is the universal central extension of $SL_d(k((t)))$,
the top row of the diagram is exact, and the homomorphisms $\rho$, $\rho_{\lam}$
and $\eps_{\lam}$ are surjective ($\rho$ is composition of the natural
map from $SL_d(k((t)))$ to $PSL_d(k((t)))$ and an isomorphism between $PSL_d(k((t)))$ and $\Ghat(A)$).
The map $\rho_{\lam}$ may or may not
be an isomorphism depending on $\lam$ (see \cite[Chapter~12]{G}).

\noindent\emph{Notational remark.}$\ $ If $A$ is a generalized Cartan matrix,
the notations $G(A)$, $\Ghat(A), U(A)$ etc. introduced in this section will
have the same meaning throughout the paper. The reference to $A$ will be omitted
when clear from the context. The last remark does not apply to Section~3
where $G$ stands for an arbitrary group.

\section{Tits' abstract simplicity theorem}

The following is a statement of the Tits
simplicity theorem for groups with a BN-pair (see \cite[Ch.~IV, No.~2.7]{B}).

\begin{theorem}
\label{thm:tits-simplicity}
Consider a quadruple $(G,B,N,U)$ where $G$ is a group,
$(B,N)$ is a BN-pair of $G$ whose associated
Coxeter system is irreducible, and $U \leq B$ is a
subgroup whose $G$-conjugates generate the entire group $G$.
\noindent Assume the following:
\negspace
\begin{enumerate}
\item[(a)] $U$ is normal in $B$ and $B=UH$, where $H=B\cap N$.
\item[(b)] $[G,G]=G$.
\item[(c)] If $\Lambda$ is a proper normal subgroup of $U$, then
$[U/\Lambda, \ U/\Lambda]
    \ \neq\
    U/\Lambda.$
\end{enumerate}
\negspace
Let $Z = \bigcap_{g\in G} gBg^{-1}$. Then the group $G/Z$ is abstractly simple.
\end{theorem}

We remark that Tits' theorem has the following ``topological" version
whose proof is identical to the ``abstract" version.

\begin{theorem}
\label{thm:tits-simplicity2}
Let $G,B,N,Z$  be as above. Assume that
$G$ is a topological group, and $B$ is a closed subgroup of $G$.
Let $U$ be a closed subgroup of $B$, and assume that
$G$ is topologically generated by the conjugates of $U$ in $G$.
Assume condition (a) above and replace (b) and (c) by conditions (b') and (c') below:

(b') $[G,G]$ is dense in $G$.

(c') If $\Lambda$ is a proper normal \bf{closed }\it subgroup of $U$, then
$[U/\Lambda, \ U/\Lambda]
    \ \neq\
    U/\Lambda.$

\noindent
Then $G/Z$ is topologically simple.
\end{theorem}

\section{Simplicity of complete Kac-Moody groups via Tits' theorem}

In this section $G$ will denote an incomplete Kac-Moody group
constructed from an \bf{indecomposable }\rm generalized Cartan matrix $A$.
Let $B,N,U$ be as in Section 2, let $\Ghat$ be the Remy-Ronan completion of $G$,
and let $\Bhat$ (resp. $\Uhat$) be the closures of
$B$ and (resp. $U$) in $\Ghat$. We shall analyze the conclusion of Tits' theorem
applied to the quadruple $(\Ghat,\Bhat,N,\Uhat)$.

The group $Z := \bigcap_{g\in \Ghat} g\Bhat g^{-1}$ is easily seen to be trivial
(see Lemma~\ref{axioms}e) and therefore, $\Ghat$ is abstractly (resp. topologically) simple
provided the hypotheses of Theorem~\ref{thm:tits-simplicity}
(resp. Theorem~\ref{thm:tits-simplicity2}) are satisfied.

We will show that the hypotheses of Theorem~\ref{thm:tits-simplicity2} are satisfied
provided $|k|>3$, thus giving a slightly different proof of Remy's theorem
on topological simplicity of $\Ghat$ \cite[Theorem 2.A.1]{Re2}.
We will also prove that hypotheses (b) and (c) of Theorem~\ref{thm:tits-simplicity} are satisfied as long
as $\Uhat$ is (topologically) finitely generated. These results will follow from
Lemma~\ref{axioms} and Lemma~\ref{lem:ClosedSubgps2} below.
Thus we will obtain the desired sufficient
condition for abstract simplicity of $\Ghat$:

\begin{theorem}
Assume that $|k|>3$. If $\Uhat$ is topologically
finitely generated, then $\Ghat$ is abstractly simple.
\label{simprel}
\end{theorem}

As an immediate consequence of Theorem~\ref{simprel}, we deduce part (b)
of Theorem~\ref{main}:

\begin{corollary}
Assume that $|k|>3$ and any $2\times 2$ submatrix of $A$ is of finite type.
Then the complete Kac-Moody group $\Ghat(A)$ is abstractly simple.
\end{corollary}
\begin{proof}
By a theorem of Abramenko~\cite{A}, under the above assumptions on $k$
and $A$, the incomplete group $U=U(A)$ is finitely generated
(as an abstract group). Thus, $\Uhat$ is automatically topologically finitely
generated, and therefore $\Ghat$ is abstractly simple by Theorem~\ref{simprel}.
\end{proof}

Before verifying Tits' simplicity axioms for the quadruple $(\Ghat,\Bhat,N,\Uhat)$,
we obtain an auxiliary result about incomplete groups.

\begin{lemma} \label{lemma:G1=G0}
The following hold: \negspace
\begin{enumerate}
\item[(a)] The group $G$ is generated by conjugates of $U$.

\item[(b)] If $|k|>3$, then $[G,G]=G$.
\end{enumerate}
\end{lemma}

\proof  (a) We know that $G$ is generated by root subgroups $\{U_{\pm\alpha_i}\}_{i\in I}$
(recall that $\{\alpha_i\}$ are simple roots). Since $U_{\alpha_i}\subset U$ and
$w_i U_{\alpha_i} w_i^{-1}=U_{w_i\alpha_i}=U_{-\alpha_i}$,
conjugates of $U$ generate $G$.

(b) Let $i\in I$ and $u\in k$, and let $g= \chi_i(u)$.
Choose $t\in k^*$ such that $t^2\neq 1$ (this is possible since $|k|>3$),
and let $v=u/(t^2-1)$.
We have
\begin{multline*}
    \chi_{i}(u)
    =
    \chi_{i}\Big( (t^{\la\alpha_i,\alpha_i^{\vee}\ra}-1 )\, v \Big)
        =
    \chi_{i}\Big( t^{\la\alpha_i,\alpha_i^{\vee}\ra} v \Big) \chi_{i}(-v)
        =
    h_{i}(t)\chi_{i}(v) h_i(t)^{-1} \chi_i(v)^{-1}
        =
    [h_i(t), \chi_i(v)].
\end{multline*}
So, $\chi_{i}(u)\in [G,G]$, and similarly one shows that $\chi_{-i}(u)\in [G,G]$.
Therefore, $[G,G]$ contains a generating set for $G$.
\QED

Now we are ready to establish Theorem~\ref{simprel} and
the corresponding statement about topological simplicity
of $\Ghat$. By Theorems~\ref{thm:tits-simplicity} and \ref{thm:tits-simplicity2},
it suffices to prove the following two results:
\footnote{The first assertion of Lemma~\ref{axioms}(b) and Lemma~\ref{lem:ClosedSubgps2}
are not needed for the proof of Theorem~\ref{simprel}.}

\begin{lemma}
The following hold:\negspace
\begin{enumerate}
\item[(a)] $\Ghat$ is generated by conjugates of $\Uhat$.
\item[(b)] $[\Ghat,\Ghat]$ is dense in $\Ghat$. Moreover, $[\Ghat, \Ghat] = \Ghat$ if $\Uhat$ is finitely generated.
\item[(c)] $\Uhat$ is normal in $\Bhat$ and $\Bhat = \Uhat (\Bhat\cap N)$.
\item[(d)] The Coxeter system of the BN-pair $(\Bhat,N)$ is irreducible.
\item[(e)] The group $Z$ is trivial.
 \end{enumerate}
\label{axioms}
\end{lemma}

\begin{lemma}  \label{lem:ClosedSubgps2}
Let $\Lambda$ be a proper normal subgroup of $\Uhat$.
If \negspace
\begin{enumerate}
\item[(a)] $\Lambda$ is closed, or

\item[(b)] $\Uhat$ is finitely generated,
\end{enumerate}
\negspace
then
\[
	[\Uhat/\Lambda,\ \ \Uhat/\Lambda]
	\ \ \neq\ \
	\Uhat/\Lambda.
\]
\end{lemma}

The proofs of Lemmas~\ref{axioms} and~\ref{lem:ClosedSubgps2}
are based on the following properties of pro-$p$ groups.

\begin{proposition} \label{prop1}
Let $K$ be a pro-$p$ group, and let $K^*$ be the closure of $[K,K]K^p$ in $K$
(where $K^p$ is the subgroup generated by $p^{\rm th}$-powers of elements of
$K$).
\negspace
\begin{enumerate}
\item[(a)] A subset $X$ of $K$ generates $K$ (topologically) if and only
if the image of $X$ in $K/K^*$ generates $K/K^*$.

\item[(b)] Suppose that $K$ is finitely generated.
Then any two minimal generating sets of $K$ have the
same cardinality.
\end{enumerate}
\negspace
\end{proposition}
\proof
(a) is proved in \cite[Proposition~25]{Se}; it also follows from
\cite[Proposition~1.9(iii)]{DSMS} and \cite[Proposition~1.13]{DSMS}.
To prove (b) note that if $K$ is finitely generated,
then $K/K^*$ can be viewed as a finite-dimensional space
over $\mathbb F_p$; let $d$ be the dimension of this space.
Let $X$ be a generating set of $K$. Clearly, $|X|\geq d$.
On the other hand, there exists a subset $Y$ of $X$,
with $|Y|=d$ such that the image of $Y$ in $K/K^*$
is a basis of $K/K^*$. By (a), $Y$ generates $K$.
Thus any minimal generating set of $K$ has cardinality $d$.
\QED

\begin{proposition}
Let $K$ be a pro-$p$ group generated by a finite set $\{a_1,\ldots, a_d\}$.
Then any element of $[K,K]$ can be written in the form
$[a_1,g_1][a_2,g_2]\ldots[a_d,g_d]$ for some $g_1,\ldots, g_d\in K$.
In particular, $[K,K]$ is closed.
\label{pro_pcomm}
\end{proposition}
\proof
The second assertion of Proposition~\ref{pro_pcomm} is the statement
of \cite[Proposition~1.19]{DSMS}. The first assertion is established
in the course of the proof of \cite[Proposition~1.19]{DSMS}.
\QED

Now we are ready to establish Lemmas~\ref{axioms} and~\ref{lem:ClosedSubgps2}.
It will be convenient to prove Lemma~\ref{lem:ClosedSubgps2} first.

\proof[Proof of Lemma~\ref{lem:ClosedSubgps2}] Suppose that $[\Uhat/\Lambda, \Uhat/\Lambda]= \Uhat/\Lambda$
or, equivalently, that $\Uhat=\Lambda[\Uhat,\Uhat]$. Then $\Lambda$
generates $\Uhat$ modulo $[\Uhat,\Uhat]$. Since $\Uhat$ is a pro-$p$ group,
$\Lambda$ generates $\Uhat$ (topologically) by Proposition~\ref{prop1}(a).
If $\Lambda$ is closed, then $\Lambda=\Uhat$, so $\Lambda$ is not proper.
Thus we proved the desired result, assuming (a).

Now assume (b) that $\Uhat$ is finitely generated.
Let $\{a_1,a_2,\ldots, a_d\}$ be a finite generating set of
$\Uhat$ contained in $\Lambda$ (such a set exists by Proposition~\ref{prop1}(b)).
Applying Proposition \ref{pro_pcomm} with $K=\Uhat$, we see that $\Lambda\supseteq [\Uhat,\Uhat]$
(since $\Lambda$ is normal in $\Uhat$). But $\Uhat=\Lambda[\Uhat,\Uhat]$,
so we conclude that $\Uhat=\Lambda$, a contradiction.
\QED

\proof[Proof of Lemma~\ref{axioms}]  {(a)} 
Let $G_1$ be the subgroup of $\Ghat$ generated by conjugates of $\Uhat$.
By Lemma~\ref{lemma:G1=G0}, $G_1$ contains $G$, whence
$G_1$ is dense in $\Ghat$. But $G_1$ is also open (hence closed) in $\Ghat$ since $G_1\supseteq \Uhat$.
Therefore, $G_1$ must be equal to $\Ghat$.

\noindent {(b)} The density of $[\Ghat,\Ghat]$ is clear
since $[\Ghat,\Ghat]\supseteq [G,G]=G$ and $G$ is dense in $\Ghat$.
Now assuming that $\Uhat$ is finitely generated, we shall
prove that $[\Ghat,\Ghat]$ is also open in $\Ghat$.
Since $\Uhat$ is a finitely generated pro-$p$ group,
its commutator subgroup $[\Uhat,\Uhat]$ is closed
by Proposition~\ref{pro_pcomm}.
So, $\Uhat/[\Uhat,\Uhat]$ is also
a finitely generated (abelian) pro-$p$ group. On the other hand,
$\Uhat$ is generated by elements of order $p$ since each root subgroup is isomorphic
to the additive group of $k$. So, $\Uhat/[\Uhat,\Uhat]$ must be finite,
and $[\Uhat,\Uhat]$ must be open in $\Uhat$ and hence in $\Ghat$. Since $\Ghat\supset \Uhat$,
we have shown that $[\Ghat,\Ghat]$ is open in $\Ghat$.

\noindent {(c)}
We know (by construction) that the corresponding results hold for incomplete
groups, that is, $U$ is normal in $B$, $B=UH$ and $H=N\cap B$. Taking the completions of both
sides of the last two equalities, and using the fact that (the images of) $H$ and $N$
in $\Ghat$ are discrete, we get $\Bhat=\Uhat H$ and $H=N\cap \Bhat$. The normality of $\Uhat$
in $\Bhat$ is clear.

\noindent {(d)}
The Coxeter group associated to $(\Bhat,N)$ is isomorphic to $W$,
and $W$ is irreducible since the matrix $A$ is indecomposable.

\noindent {(e)}
Recall that $X^+$ denotes the building associated with the incomplete BN-pair $(B,N)$,
and let $\mathcal C$ be the chamber of $X^+$ whose stabilizer in $G$ is $B$. It follows directly
from definitions that the stabilizer of $\mathcal C$ in $\Ghat$ is $\Bhat$.
So, $Z=\bigcap_{g\in \Ghat}g\Bhat g^{-1}$ consists of elements
which stabilize all chambers in $X^+$ (recall that $G$ acts transitively on the set of chambers).
Therefore $Z=\{1\}$.
\QED

\section{``Relative" Kac-Moody groups}

Let $A=(a_{ij})_{i\in I}$ be a generalized Cartan matrix. 
As before, let $\Pi=\{\alpha_1,\ldots,\alpha_l\}$ be the set of simple roots.
Recall that given a subset $J$ of $I$, we denote by $A_J$ the $|J|\times |J|$ matrix $(a_{ij})_{i,j\in J}$.

One can associate two \emph{a priori} different groups to the matrix $A_J$: the
usual Kac-Moody group $G(A_J)$, and the ``relative" Kac-Moody group
$G_J$ which is defined as the subgroup of $G(A)$ generated
by $\{U_{\pm \alpha_i}\}_{i\in J}$. The first result of this section
asserts that these two groups are canonically isomorphic
(Proposition~\ref{relative1}). 
Next, one can consider two different topologies on $G(A_J)$ -- the usual
building topology and the topology induced from the building topology on $G(A)$ via the
above isomorphism $G(A_J)\cong G_J$.
While the two topologies may not be the same, the main result of this section is
that their restrictions to the subgroup $U(A_J)$ are the same, provided $A$ is symmetric
(Theorem~\ref{relative2}). We believe that Theorem~\ref{relative2} holds without the
assumption that $A$ is symmetric, but we are unable to prove it at the present time.
In the next section we use Theorem~\ref{relative2} to reduce the problem
of finite generation of $\Uhat(A)$ to a certain question about the groups $\Uhat(A_J)$
where $J$ runs over all subsets of cardinality $2$, once again assuming
that $A$ is symmetric.

In order to state  our results precisely, we introduce the following notation.
We set $G=G(A)$, $\Phi=\Phi(A)$, $W=W(A)$, $U=U(A)$ and $U^- = U^-(A)$. Fix a
subset $J\subseteq I$, and let $\Pi_J=\{\alpha_i\}_{i\in J}$, let $W_J$ be the
subgroup of $W$ generated by $\{w_i\}_{i\in J}$
and $\Phi_J=W_J (\Pi_J)$. Also let $\Phi_J^{\pm}=\Phi_J\cap \Phi^{\pm}$.
Clearly, $W_J$ (resp. $\Phi_J$) can be canonically identified
with $W(A_J)$ (resp. $\Phi(A_J)$).

We have already defined the group $G_J$. Note that by relations (R3),
$G_J\supset U_{\alpha}$ for every $\alpha\in \Phi_J$.
Also introduce the subgroups
\[
	U_J:=\{U_{\alpha}\mid \alpha\in \Phi_J^+\}
	\ \text{ and } \
	U_J^-:=\{U_{\alpha}\mid \alpha\in \Phi_J^-\}
\]
of $G_J$. We will see shortly that
$U_J=G_J\cap U$ and $U_J^-=G_J\cap U^-$.

\begin{proposition}
There is a natural isomorphism $\phi: G(A_J)\to G_J$.
Moreover, $\phi$ maps $U(A_J)$ onto $U_J$ and $B(A_J)$ onto $B_J$.
\label{relative1}
\end{proposition}

In view of this proposition, we can identify $U_J$ (resp. $G_J$) with $U(A_J)$ (resp. $G(A_J)$).
Let $\Uhat_J$ be the closure of $U(A_J)$ in $\Ghat(A_J)$, as before,
and let $\Ubar_J$ be the closure of $U_J$ in $\Ghat$. 

\begin{theorem} \label{relative2}
Assume that $A$ is symmetric. The groups $\Uhat_J$ and $\Ubar_J$ are (topologically) isomorphic.
\end{theorem}

\proof[Proof of Proposition~\ref{relative1}]
As before, identify $\Phi(A_J)$ with the subset $\Phi_J$ of $\Phi$. To distinguish between
generators of $G(A_J)$ and $G(A)$, we use the symbols
$\{x_{\alpha}(u)\mid \alpha\in \Phi(A_J), u\in k\}$ for the generators
of $G(A_J)$ (the generators of $G(A)$ are denoted $\{\chi_{\alpha}(u) \}$ as usual).
From the defining presentation of Kac-Moody groups,
it is clear that there exists a map $\phi:G(A_J)\to G$
such that $\phi(x_{\alpha}(u))=\chi_{\alpha}(u)$ for $\alpha\in \Phi(A_J), u\in k$.
Clearly, $\phi(G(A_J))=G_J$, $\phi(U(A_J))=U_J$ and $\phi(B(A_J))=B_J$,
so we only need to show that $\phi$ is injective. We proceed in several steps.

\it{Step 1: $\phi$ is injective on $U(A_J)$. }\rm
Let $\psi: U(A)\to U(A_J)$ be the unique homomorphism such that
$$\psi (\chi_{\alpha}(u))=
\left\{
\begin{array}{ll}
x_{\alpha}(u)& \mbox{ if } \alpha\in \Phi_J,\\
1& \mbox{ if } \alpha \not\in\Phi_J.\\
\end{array}
\right.
$$
The fact that such a homomorphism exists follows immediately from Theorem~\ref{pres_U}.
It is also clear that $\psi\phi(g)=g$ for any $g\in U(A_J)$, whence
the restriction of $\phi$ to $U(A_J)$ must be injective.

\it{Step 2: $\phi$ is injective on $H(A_J)$. }\rm
This follows directly from the fact that relations (R6)-(R7) are defining
relations for the tori $H(A_J)$ and $H(A)$.

\it{Step 3: $\phi^{-1}(B_J)=B(A_J)$ where $B_J=B\cap G_J$. }\rm
It is clear that $\phi^{-1}(B_J)\supseteq B(A_J)$.
Suppose that $\phi^{-1}(B_J)$ is strictly larger than $B(A_J)$.
Since $(B(A_J), N(A_J))$ is a BN-pair of $G(A_J)$, we conclude that $\phi^{-1}(B_J)$
is of the form $B(A_J) W(A_K) B(A_J)$ for some non-empty subset $K\subseteq J$.
This would mean that $B_J=\phi(B(A_J) W(A_K) B(A_J))$ contains at least one of the generators
of $W$, which is impossible since $B\cap N=H$.

\it{Step 4: Conclusion. }\rm
Let $K=\Ker\phi$. Note that $K\subset B(A_J)$ since $B(A_J)$ is the full preimage
of $B_J$ under $\phi$. Take any $g\in K$ and write it as $g=uh$ where $u\in U(A_J)$ and $h\in H(A_J)$.
Then $\phi(u)=\phi(h)^{-1}$. On the other hand, it is clear that $\phi(u)\in U$ and $\phi(h)\in H$.
Since $H\cap U={1}$ and $\phi$ is injective on both $U(A_J)$ and $H(A_J)$,
we conclude that $u=h=1$. Therefore, $K$ is trivial.
\QED

Before proving Theorem~\ref{relative2}, we state two lemmas, 
which will be established at the end of the section.

\begin{lemma}
Let $\Psi_J=\{\beta\in\Phi^+: \la\beta, \alpha_j^{\vee}\ra\leq 0\mbox{ for any }j\in J\}$.
The following hold:

(a) For any $\gamma\in \Phi^+\backslash \Phi_J^+$ there exist $w\in W_J$ and $\beta\in \Psi_J$
such that $\gamma=w\beta$.

(b) Assume that $A$ is symmetric. Then for any $\gamma\in \Psi_J $ and $\beta\in \Phi_J^+$, the root groups
$U_{-\gamma}$ and $U_{\beta}$ commute (elementwise).
\label{revised}
\end{lemma}

\begin{lemma}
Let $C$ be a generalized Cartan matrix. Then the building topology on
$U(C)$ is given by the subbase $\{gU(C) g^{-1}\cap U(C)\}_{g\in G(C)}$.
\label{T2=T2'}
\end{lemma}

\proof[Proof of Theorem~\ref{relative2}]
By definition, $\Ubar_J$ and $\Uhat_J$ are the completions of $U_J$
with respect to the topologies $(\T_1)$ and $(\T_2)$, respectively, where
$(\T_1)$ is given by the subbase $\{gBg^{-1}\cap U_J\}_{g\in G}$ and
$(\T_2)$ is given by the subbase $\{gB_J g^{-1}\cap U_J\}_{g\in G_J}$. We have to show that
$(\T_1)$ and $(\T_2)$ coincide.

The inequality $(\T_1)\geq (\T_2)$ is clear.
Indeed, for any $g\in G_J$ we have $gBg^{-1}\cap\, U_J=gB_J g^{-1}\cap\, U_J$
since $gB_J g^{-1}=g(B\cap G_J)g^{-1}=gBg^{-1}\cap G_J$.

Now we prove the reverse inequality.
In view of the natural isomorphism $U(A_J)\cong U_J$,
Lemma~\ref{T2=T2'} applied with $C=A_J$ reduces
the proof of the inequality $(\T_2)\geq (\T_1)$ to
the following statement:

\begin{claim}
Given $g\in G$, there exists a finite set $T\subset G_J$
such that
\[
    gBg^{-1} \cap U_J\quad\  \supseteq \quad\
    \bigcap_{t\in T} t\, U_J\, t^{-1}\ \cap\ U_J.
\]
\label{***}
\end{claim}
\vskip -.5cm
Fix $g\in G$. By the Birkhoff decomposition, $g=g_- w g_+$
for some $g_+\in B$, $w\in W$ and $g_-\in U^-$.
We will show that there exist $g_1, \ldots, g_{k}\in G_J$ and $v\in W_J$
such that
\negspace
\begin{enumerate}
\item[(a)] if $x\in U_J$ is such that $g_i^{-1} x g_i\in U_J$ for $1\leq i\leq k$,
then $g_-^{-1} x g_- = g_k^{-1} x g_k$ (in particular, $g_-^{-1} x g_-\in U_J$),

\item[(b)] if $y\in U_J$ is such that $v^{-1}yv\in U_J$, then $w^{-1}yw\in U$.
\end{enumerate}
\negspace

First, let us see why (a) and (b) will imply Claim~\ref{***}. Indeed, let $\{g_i\}$ and $v$
be as above, and set $T=\{g_i\}_{i=1}^k\cup \{g_k v\}$. Let $x\in \bigcap_{t\in T} t\, U_J\, t^{-1}\cap U_J$.
Then by (a), $g_-^{-1} x g_- =g_k^{-1} x g_k$. Now let $y=g_-^{-1} x g_-$. Then
$v^{-1}yv=(g_k v)^{-1}x (g_k v)\in U_J$ by the choice of $T$, so applying (b) we get
that $w^{-1} g_-^{-1}x g_- w\in U$. Finally, $g^{-1}xg=g_+^{-1} w^{-1} g_-^{-1}x g_- w g_+\in U$
since $g_+\in B$ and $B$ normalizes $U$. So, $x\in g U g^{-1}$.

\vskip .1cm
\emph{Proof of (a).}
By definition of $U^-$, we can write $g_-$ in the form $t_1\ldots t_s$ such that
each $t_i$ lies in $U_{-\gamma_{i}}$ for some $\gamma_i\in \Phi^+$. By Lemma~\ref{revised}(a), for each
$i$ such that $\gamma_i\not\in \Phi_J$, there exists $v_i\in W_J$ and $\beta'_i\in \Psi_J$ such that 
$\gamma_i=v_i\beta'_i$. Then $t_i\in U_{-\gamma_{i}}=U_{v_i(-\beta'_i)}=v_i U_{-\beta'_i}v_i^{-1}$.
Thus, we can write $g_- =x_1\ldots x_k$ (with $k\leq 3s$), where for each $1\leq i\leq k$ either
$x_i\in U_{-\beta_i}$ with $\beta_i\in \Psi_J$ or $\beta_i\in \Phi_J$, or
$x_i\in W_J$. Note that for each $1\leq i\leq k$ either $x_i\in U_{-\beta_i}$ with $\beta_i\in \Psi_J$,
or $x_i\in G_J$.

Now define
$y_0,\ldots, y_k\in G$
and $g_0,\ldots, g_k\in G_J$ inductively: $y_0=1$; $y_i=y_{i-1}x_i$
for $1\leq i\leq k$; $g_0=1$, and for $1\leq i\leq k$ we set
$g_i=g_{i-1}x_i$ if $x_i\in G_J$ and $g_i=g_{i-1}$ if
$x_i\not\in G_J$. Note that $y_k=g_-$ and $g_i\in G_J$ for each $i$.

Suppose that $x$ satisfies the hypotheses of (a).
We shall prove by induction that \mbox{$y_i^{-1}x y_i=g_i^{-1}x g_i$}
for all $i$. The base case $i=0$ is trivial. Now let $i>0$, and suppose that
\mbox{$y_{i-1}^{-1}x y_{i-1}=g_{i-1}^{-1}x g_{i-1}$}.
If $x_i\in G_J$, then by construction $y_{i-1}^{-1}y_i=g_{i-1}^{-1}g_i=x_i$, so the induction step is clear.

If $x_i\not\in G_J$, then by construction $x_i\in U_{-\beta_i}$ for some $\beta_i\in\Psi_J$. 
By Lemma~\ref{revised}(b), $U_{-\beta_i}$ centralizes $U_J$, 
since $U_J$ is generated by $\{U_{\gamma}\}_{\gamma\in\Phi_J^+}$.
By the induction hypotheses, $y_{i-1}^{-1}x y_{i-1}=g_{i-1}^{-1}x g_{i-1}\in U_J$.
Therefore, $$y_i^{-1}x y_i=x_i^{-1}(y_{i-1}^{-1}x y_{i-1})x_i=y_{i-1}^{-1}x y_{i-1}=
g_{i-1}^{-1}x g_{i-1}=g_i^{-1}x g_i.$$ This completes the induction step and hence the proof of (a).

\vskip .1cm
\emph{Proof of (b).} The following argument was suggested to us (in a slightly
different form)
by Bertrand Remy and Pierre-Emmanuel Caprace.
First of all, the statement of (b) is clearly implied by the following: for any $w\in W$
there exists $v\in W_J$ such that $wv U_J (wv)^{-1}\subseteq U$.
We know that $U_J$ is generated by the root subgroups
$\{U_{\alpha}\}_{\alpha\in\Phi_J^+}$.
Since $z U_{\alpha}z^{-1}=U_{z\alpha}$ for every $z\in W$, it suffices to prove
the following statement about root systems:
\vskip .2cm
\centerline{For any $w\in W$
there exists $v\in W_J$ such that $wv(\alpha_i)>0$ for all $i\in J$.}

Recall that the simple roots $\{\alpha_i\}$ are linearly independent elements
of the $\mathbb Q$-vector space $\mathfrak h^*$. Let
$\mathfrak h^*_{\mathbb R} = \mathfrak h^*\otimes_{\mathbb Q} \mathbb R$.
Let $R=\{\beta\in \mathfrak h^*_{\mathbb R}\mid \la\beta,\alpha_i^{\vee}\ra> 0 \mbox{ for all } i\}$
be the fundamental chamber and $R_J=
\{\beta\in \mathfrak h^*_{\mathbb R}\mid \la\beta,\alpha_i^{\vee}\ra> 0 \mbox{ for all } i\in J\}$
be the $J$-relative fundamental chamber (in the sense of \cite[5.7]{MP}). Now pick any $x\in R$
(it is well known that $R$ is non-empty).
By \cite[5.7, Proposition 5(i)]{MP}, the union of $W$-translates of $R$ is contained in the union of
$W_J$-translates of $R_J$. Therefore, given $w\in W$, there exists $v\in W_J$
such that $v^{-1}w^{-1}x\in R_J$. For any $i\in J$ we have $\la v^{-1}w^{-1}x, \alpha_i^{\vee}\ra>0$,
whence $\la x, wv\alpha_i^{\vee} \ra>0$. Since $x\in R$, the last inequality implies that
$wv\alpha_i^{\vee}$ must be a positive coroot, whence $wv\alpha_i$ is a positive root.
The proof is complete.
\QED
\begin{proof}[Proof of Lemma~\ref{revised}]
(a) This result is probably well known and follows easily from \cite[5.7, Proposition 5]{MP},
but since a direct proof is very short, we present it here.
The proof is by induction on $\height(\gamma)$.

If $\gamma\in\Phi^+\backslash \Phi_J^+$ and $\height(\gamma)=1$, then $\gamma=\alpha_i$ for some $i\not\in J$,
whence for any $j\in J$ we have $\la \gamma,\alpha_j^{\vee}\ra=a_{ji}\leq 0$, so $\gamma\in \Psi_J$.

Now take any $\gamma\in\Phi^+\backslash \Phi_J^+$ and assume that Lemma~\ref{revised}(a) holds for
any root $\gamma'\in\Phi^+\backslash \Phi_J^+$ with $\height(\gamma')<\height(\gamma)$.
If $\la \gamma,\alpha_j^{\vee}\ra\leq 0$ for any $j\in J$, then $\gamma\in\Psi_J$.
If $\la \gamma,\alpha_j^{\vee}\ra>0$ for some $j\in J$, let $\gamma'=w_j\gamma$.
Then $\gamma'=\gamma-\la \gamma,\alpha_j^{\vee}\ra\alpha_j$ is a root
of smaller height than $\gamma$, and $\gamma'\in \Phi^+\backslash \Phi_J^+$ since the set
$\Phi^+\backslash \Phi_J^+$ is $W_J$-invariant. Thus, there exist $w\in W_J$ and $\beta\in \Psi_J$
such that $\gamma'=w\beta$, whence $\gamma=w_j \gamma'=(w_j w)\beta$.
\vskip .12cm

(b) Since $A$ is symmetric, there exists a symmetric $W$-invariant
bilinear form $(\cdot,\cdot)$ on $\mathfrak h^*$ such that 
$(\alpha,\beta)=\la\alpha,\beta^{\vee}\ra$ for any $\alpha,\beta\in\Phi$;
in particular, $(\alpha,\alpha)=2$  for any $\alpha\in \Phi$ (see \cite[Chapter~2]{Kac}).

Now let $\gamma\in \Psi_J $ and $\beta\in \Phi_J^+$. Then 
the set $\dbZ_{>0}(\beta)+\dbZ_{>0}(-\gamma)$ does not contain any
real roots. Indeed, $(\gamma,\beta)\leq 0$ by definition of $\Psi_J$,
so for any $i,j>0$ we have
$(i\beta-j\gamma,i\beta-j\gamma)=2(i^2+j^2-ij(\gamma,\beta))\geq 2(i^2+j^2)\geq 4$.
Thus the pair $\{\beta,-\gamma\}$ is prenilpotent and moreover, the corresponding root groups
commute by relations (R2).
\end{proof}

\proof[Proof of Lemma~\ref{T2=T2'}]
It will suffice to show that
$gU(C) g^{-1}\cap U(C)=gB(C) g^{-1}\cap U(C)$ for any $g\in G(C)$.
Thus we fix $g\in G(C)$ and $x\in U(C)$ such that
$g^{-1}xg\in B(C)$. Then $g^{-1}xg=hu$ where $u\in U(C)$ and $h\in H(C)$.
We need to show that $h=1$.

Assume that $h\neq 1$. Then we can choose a matrix $D$ such that $C$ is a submatrix of $D$ and $h$
does not lie in the center of $G(D)$ (as before, we identify $G(C)$ with a subgroup of $G(D)$).
The existence of such a $D$ follows easily from defining relations (R4) -- we just have to ensure
that $h$ acts non-trivially by conjugation on one of the root subgroups.

Now let $(\T)$ be the building topology on $G(D)$.
Let $\Ubar(C)$ (resp. $\Ghat(D)$, $\Uhat(D)$) be the completion of $U(C)$ (resp. $G(D)$, $U(D)$)
with respect to $(\T)$. By Theorem~\ref{collection}b), $\Uhat(D)$
is a pro-$p$ group, hence so is $\Ubar(C)$.
Therefore, $x^{p^n}\to 1$ in $(\T)$ as $n\to \infty$, whence $(g^{-1}xg)^{p^n}\to 1$ in $(\T)$
as well. On the other hand, $(g^{-1}xg)^{p^n}=(hu)^{p^n}=h^{p^n}u_n$  where $u_n\in U(C)$.
Since the order of $h$ is finite, prime to $p$ and $\Ubar(C)$ is compact, there exists a subsequence
$\{n_k\}$ such that $h^{p^{n_k}}=h$ for all $k$
and $\{u_{n_k}\}$ converges to some element $v\in \Ubar(C)$. Thus $hv=1$ in $\Ghat(D)$. Since
$\Ubar(C)$ is pro-$p$, we conclude that $h=1$ in $\Ghat(D)$.
So, $h$ lies in the kernel of the natural map $G(D)\to \Ghat(D)$ which, as we know,
coincides with the center of $G(D)$. This contradicts our initial assumption.
\QED

\section{Finite generation of $\Uhat$: reduction to the rank 2 case}
In this section we prove that finite generation of $\Uhat$
is essentially determined by rank $2$ subsystems of $\Phi$, provided $A$ is symmetric.
\begin{definition}
Let $A$ be a generalized Cartan matrix and $\Uhat=\Uhat(A)$. We say that $\Uhat$ is 
\it{well behaved }\rm
if for any non-simple root $\gamma\in \Phi(A)^+$ we have $U_{\gamma}\subseteq [\Uhat, \Uhat]$.
\end{definition}

By Proposition~\ref{prop1}(a), if $\Uhat$ is well behaved then $\Uhat$ is topologically generated
by simple root subgroups $\{U_{\alpha}\}_{\alpha\in \Pi(A)}$ (in particular, $\Uhat$
is topologically finitely generated).

For the rest of this section we fix a matrix $A=(a_{ij})_{i,j\in I}$, and we write $\Uhat=\Uhat(A)$
and $\Phi=\Phi(A)$.
For a subset $J$ of $I$ we define $\Phi_J$ as in the previous section.
Recall that $\{\alpha_i\}_{i\in I}$ are simple roots, and $w_i$ is the
reflection associated with $\alpha_i$.

\begin{theorem} Suppose that $A$ is symmetric, and
for any subset $J\subseteq I$ of cardinality $2$,
the group $\Uhat(A_J)$ is well behaved.
Then $\Uhat$ is also well behaved and hence (topologically)
finitely generated.
\label{fingen}
\end{theorem}

Theorem~\ref{fingen} is an easy consequence of the following lemma.
\begin{lemma}
Let $\gamma\in \Phi^+$ be a non-simple root.
There exist simple roots $\alpha_i,\alpha_j$ and $w\in W$
such that $w\alpha_i>0$, $w\alpha_j>0$ and
$\gamma=w\alpha$ for some non-simple root
$\alpha\in\Phi^+_{\{i,j\}}$.
\label{prenil2}
\end{lemma}

\proof [Proof of Lemma~\ref{prenil2}]

If $\gamma$ lies in a subsystem generated by two simple roots,
that is, $\gamma=n\alpha_i+m\alpha_j$ for some $n,m\in\dbN$ and $i,j\in I$,
the assertion is obvious (we can take $w=1$).
From now on assume that
$\gamma=\sum_{i=1}^l n_i\alpha_i$ where at least three $n_i$'s are nonzero.
We will prove the lemma by induction on $\height(\gamma)$.

We note that there exists $k\in I$ such that
$\height(w_{k}\gamma)<\height(\gamma)$. Indeed, if $\height(w_{i}\gamma)\geq \height(\gamma)$
for all $i$, then $\la \gamma,\alpha_i^{\vee}\ra\leq 0$ for all $i$
and hence $\la \alpha_i,\gamma^{\vee}\ra\leq 0$ for all $i$.
But $\gamma$ is a linear combination
of $\alpha_i$ with nonnegative coefficients, so we must have $\la\gamma,\gamma^{\vee}\ra\leq 0$.
The latter is impossible since $\la\gamma,\gamma^{\vee}\ra=2$.

Since $\gamma$ does not lie in a rank two subsystem,
$w_k\gamma$ is not simple. By induction,
there exist simple roots $\alpha_i$ and $\alpha_j$, $w\in W$
and $\alpha\in\Phi_{\{i,j\}}^+$ such that
$w\alpha_i>0$, $w\alpha_j>0$ and $w_k\gamma=w\alpha$.

Note that $\gamma=w_k w\alpha$.
If both $w_k w\alpha_i$ and $w_k w\alpha_j$ are positive,
we are done. Suppose, this is not the case. Then we must have
$w\alpha_i=\alpha_k$ or $w\alpha_j=\alpha_k$, and without loss of
generality we assume that $w\alpha_i=\alpha_k$.
We have
\begin{multline*}
\gamma=w_k (w\alpha)=w\alpha - \la w\alpha,\alpha_k^{\vee}\ra\alpha_k=
w\alpha - \la w\alpha,(w\alpha_i)^{\vee}\ra w\alpha_i=\\
w\alpha - \la w\alpha,w(\alpha_i^{\vee})\ra w\alpha_i=
w\alpha-\la \alpha,\alpha_i^{\vee}\ra w\alpha_i=
w(\alpha-\la \alpha,\alpha_i^{\vee}\ra \alpha_i)=
w(w_i\alpha).
\end{multline*}

Since $\alpha\in\Phi^+_{\{i,j\}}$ and $\alpha$ is not simple,
we have $w_i\alpha\in \Phi^+_{\{i,j\}}$. The proof will
be complete if we show that $w_i\alpha$ is not simple.
If $w_i\alpha$ is simple, then $\height(w_i\alpha)<\height(\alpha)$,
so we must have $\la\alpha,\alpha_i^{\vee}\ra>0$. Therefore,
$\la w\alpha,\alpha_k^{\vee} \ra=\la w\alpha,w\alpha_i^{\vee} \ra>0$, whence
$$\height(\gamma)=\height(w_k w\alpha)<\height(w\alpha)=\height(w_k\gamma),$$ contrary
to our assumptions.
\QED

\proof [Proof of Theorem~\ref{fingen}]
Fix a non-simple positive real root $\gamma$. We have to show
that $U_{\gamma}\in [\Uhat,\Uhat]$.
By Lemma~\ref{prenil2},
there exist $\alpha_i,\alpha_j\in\Pi$, $w\in W$
and a non-simple root $\alpha\in\Phi^+_{\{i,j\}}$
such that $w\alpha_i>0$, $w\alpha_j>0$ and
$\gamma=w\alpha$. Let $J=\{i,j\}$ and $V=\Uhat_{J}$.
Since $V\cong \Uhat(A_J)$ by Theorem~\ref{relative2},
the hypotheses of Theorem~\ref{fingen} imply that
$V$ is topologically generated by
$U_{\alpha_i}\cup U_{\alpha_j}$.
Since $w U_{\beta} w^{-1}=U_{w\beta}$ for any $\beta\in\Phi$,
we conclude that $ w V w^{-1}\subseteq \Uhat$. By hypotheses, we also know that
$U_{\alpha}\subset [V,V]$. Therefore,
$U_{\gamma}=U_{w\alpha}=w U_{\alpha} w^{-1}\subset
[w V w^{-1}, w V w^{-1}]\subseteq [\Uhat,\Uhat]$.
\QED

In the next section we establish a sufficient condition
for finite generation of $\Uhat$ in the rank two case:

\begin{theorem}
Let $C$ be a $2\times 2$ generalized Cartan matrix.
Assume that either $C$ is finite and $|k|>3$, or that $C$ is affine,
$|k|>3$ and $p>2$.
Then the group $\Uhat(C)$ is well behaved.
\label{rank2}
\end{theorem}

Before proving Theorem~\ref{rank2} we explain how to deduce
Theorem~\ref{main}(a).
Theorem~\ref{mainfg} below is a direct consequence of Theorem~\ref{fingen}
and Theorem~\ref{rank2}. Theorem~\ref{main}(a) is obtained by combining
Theorem~\ref{mainfg} and Theorem~\ref{simprel}.

\begin{theorem}\label{mainfg}
Let $A$ be an indecomposable generalized Cartan matrix. Suppose that 
$A$ is symmetric and and any $2\times 2$ submatrix of $A$ is of finite or affine type,
$|k|>3$ and $p=\char(k)>2$.
Then the group $\Uhat(A)$ is finitely generated.
\end{theorem}

\begin{remark} As one can see from the proofs of Theorem~\ref{fingen} and Theorem~\ref{mainfg}, 
the assumption `$A$ is symmetric' was only needed to apply Theorem~\ref{relative2}.
We believe that Theorem~\ref{relative2} can be proved (by a completely different method)
without the assumption `$A$ is symmetric', but under the assumption `$A_J$ is of affine type',
using explicit realization of affine Kac-Moody groups and solution to the congruence subgroup
problem. If the latter is achieved, one would be able to eliminate the assumption `$A$ is symmetric'
from the statements of Theorems~\ref{fingen} and \ref{mainfg} and Theorem~\ref{main}(a). 
\end{remark}

\section{Finite generation: the rank 2 case}
\label{sec:ranktwo}

In this section we prove Theorem~\ref{rank2}. In order to finish the proof of
Theorem~\ref{main}(a), it would have been enough to prove Theorem~\ref{rank2} for symmetric matrices of finite or
affine type. However, we decided to prove the full version of Theorem~\ref{rank2},
including the computationally demanding case of non-symmetric affine matrices, because
it may have applications to other problems, including possible generalization
of Theorem~\ref{main}(a) (as explained at the end of Section~6).
\vskip .1cm

If $C$ is of finite type, the group $G(C)$ is finite, so
$\Uhat(C)$ is isomorphic to $U(C)$. Theorem~\ref{rank2} in this case
follows easily from defining relations of type (R2)
(see \cite{Mo}, where the coefficients $\{C_{mn\alpha\beta}\}$
are computed explicitly in terms of the root system generated by $\alpha$ and $\beta$).
Before considering the affine case, we make
some general remarks about affine Kac-Moody groups.

A good general reference for incomplete affine Kac-Moody groups
is the thesis of Ramagge \cite{Ra1} (see also \cite{Ra2}, \cite{Ra3}).
In particular, using results of \cite{Ra1}, one can obtain
an explicit realization of incomplete twisted affine Kac-Moody groups.
An explicit realization of complete twisted affine Kac-Moody groups
is probably known, however we are unaware of a proof in the literature.
For completeness, we shall demonstrate such realization
in the case of $2\times 2$ matrices - see Proposition~\ref{app}.

Let $C$ be an affine matrix (of arbitrary size). Then the incomplete group $G(C)$
modulo its finite center is isomorphic to the group of fixed points
of a finite order (possibly trivial) automorphism $\omega_C$ of the group
of $k[t,t^{-1}]$-points of some simply-connected Chevalley group.
The automorphism $\omega_C$ is non-trivial if and only if $C$ is a
\emph{twisted} affine matrix.\footnote{For a classification of twisted affine
matrices see \cite[Chapter 4]{Kac}. } The complete group
$\Ghat(C)$ has an analogous
description where the ring $k[t,t^{-1}]$ is replaced by the field $k((t))$.
Furthermore, $\Ghat(C)$ is isomorphic to the group of $k((t))$-points of some simple
algebraic group $\mathbb G_C$ defined over $k((t))$: if $C$ is non-twisted, $\mathbb G_C$
is the Chevalley group mentioned above; if $C$ is twisted, $\mathbb G_C$ is non-split.

One may ask if every simple algebraic group over $k((t))$ is isogenous
to one of the form $\Ghat(C)$ for some $C$. The answer is `no' since the
groups $\Ghat(C)$ are always residually split. Moreover, there is
a bijective correspondence between isogeny classes of residually
split simple algebraic groups over $k((t))$ and groups of the form
$\Ghat(C)$, with $C$ affine.
If $C$ is an affine matrix of type $X_n^{(r)}$ (in the notation of \cite[Chapter~4]{Kac}),
then $\Ghat(C)$ is a residually split group whose isogeny class is
given by Tits' index of the form ${}^{r}\!X_{n,m}^{(d)}$ for some $m,d$ in the
notation of \cite{Ti} (these conditions determine the isogeny class uniquely).

Now assume that $C$ is a $2\times 2$ affine matrix. Up to isomorphism of the corresponding Kac-Moody groups,
there are only two possibilities:
\[
C=
\begin{pmatrix}
2 &-2\\
-2& 2\\
\end{pmatrix}
\ \ \text{(non-twisted case)}, \quad \text{ and }
\quad C=
\begin{pmatrix}
2 &-4\\
-1& 2\\
\end{pmatrix}
\ \ \text{ (twisted case) }.
\]

In the first case $C$ has type $A_1^{(1)}$, the Tits index is
${}^{1}\!A_{1,1}^{(1)}$,
and $\Ghat(C)$ is isomorphic to $PSL_2(k((t)))$. In the second case
$C$ has type $A_2^{(2)}$, the Tits index is ${}^{2}\!A_{2,1}^{(1)}$, and
$\Ghat(C)$ is isomorphic to $SU_3(k((t)),h,\sigma)$, where $\sigma$ is
a $k$-automorphism of $k((t))$ of order $2$ and $h$ is
a hermitian form in three variables relative to $\sigma$ (see \cite{Ti}).
If $p\neq 2$, we can (and will) assume that
\[
	\sigma(t)= -t
	\ \text{ and } \
	h((x_1,x_2,x_3), (y_1,y_2,y_3))=x_1\sig(y_3)-x_2\sig(y_2)+x_3\sig(y_1) .
\]
Then
\[
	SU_3(k((t)),h,\sigma)=\{g\in SL_3(k((t)))\mid J^{-1}g^* J=g^{-1}\}
	\ \text{ where } \
	J=
	\left(
		\begin{smallmatrix} 0&0&1\\  0&-1&0\\  1&0&0\\ \end{smallmatrix}
	\right)
\]
and $g\mapsto g^*$ is the map from $SL_3(k((t)))$ to $SL_3(k((t)))$
obtained by applying $\sigma$ to each entry of $g$ followed by transposition.

The proof of Theorem~\ref{rank2} in both cases follows the same strategy,
but the twisted case requires more computations.
Our starting point is the following obvious lemma:

\begin{lemma} \label{wellbeh}
Let $P$ be a profinite group, and let
$P_1\supset P_2\supset ...$ be a descending chain of open normal subgroups of $P$
which form a base of neighborhoods of identity. For each $i\in\dbN$
choose elements $g_{i,1},\ldots, g_{i,n_i}\in P_i$ which generate $P_i$ modulo $P_{i+1}$.
Let $K$ be a closed subgroup of $P$ such that $g_{i,j}\in KP_{i+1}$ for all $i$ and $j$.
Then $K$ contains $P_1$.
\QED\end{lemma}

We shall apply this lemma with $P=\Uhat$, $K=[P,P]$ and a certain filtration
$\{P_i\}$ of $P$ satisfying the above conditions and such that $P_1$ contains
$U_{\gamma}$ for every non-simple positive root $\gamma$. Clearly, this will prove that
$\Uhat$ is well behaved, so we only need to show the existence of a filtration with
required properties.

\emph{Case 1: }
$C=
\begin{pmatrix}
2 &-2\\
-2& 2\\
\end{pmatrix}$.

In this case,
$\Phi^+=\{(n+1)\alpha_1+n\alpha_2, n\alpha_1+(n+1)\alpha_2\}_{n \in\dbZ_{\geq 0}}$.
The group $\Uhat$ can be embedded
into $SL_2(k[[t]])$ via the following map:
$$\chi_{(n+1)\alpha_1+n\alpha_2}(u)\mapsto
e_n(u):=
\left(
\begin{array}{ll}
1 &u t^{n}\\
0& 1\\
\end{array}
\right),\quad
\chi_{n\alpha_1+(n+1)\alpha_2}(u)\mapsto
f_{n+1}(u):=
\left(
\begin{array}{ll}
1 &0\\
u t^{n+1}& 1\\
\end{array}
\right),
$$ for all $n\in\dbZ_{\geq 0}$ and $u\in k$.
For each $n\geq 0$ let $E_n=\{e_n(u)\mid u\in k\}$ and $F_n=\{f_n(u)\mid u\in k\}$.
We also define the elements $\{h_n(u)\mid n\geq 1,u\in k \}$ by
\[
	h_n(u)
	=
	\begin{pmatrix} 1+u t^n &0\\ 0& (1+u t^n)^{-1}\\ \end{pmatrix} .
\]

Under the above identification, $\Uhat$ consists of matrices in $SL_2(k[[t]])$
whose reduction$\mod t$ is upper-unitriangular.
Now define the filtration $\{P_n\}$ as follows:
$P_1=P_2 E_1$, and for $n\geq 2$ set
\[
	P_n = SL_2^n (k[[t]]) = \{g\in SL_2(k[[t]])\mid g\equiv 1\mod t^n\}.
\]
Clearly, $P_1$ contains all non-simple positive root subgroups.

For each $n\geq 2$, $P_n$ is generated modulo $P_{n+1}$
by the elements $\{e_{n}(u), f_n(u), h_n(u)\mid u\in k\}$,
and
$P_1$ is generated modulo $P_2$ by $\{e_1(u)\}$.
Direct computation shows that
$e_{n}(u)\equiv [h_1(1),e_{n-1}(u/2)]\mod P_{n+1}$ for $n\geq 1$,
$f_n(u)\equiv [h_1(1),f_{n-1}(-u/2)]\mod P_{n+1}$ for $n\geq 2$ and
$h_n(u)\equiv [e_1(1), f_{n-1}(u)]\mod P_{n+1}$ for $n\geq 2$.
So, all the hypotheses of Lemma~\ref{wellbeh} are satisfied,
and we are done with Case 1.

\emph{Case 2: }
$C=
\begin{pmatrix}
2 &-4\\
-1& 2\\
\end{pmatrix}$.

\noindent
Let $\delta=2\alpha_1+\alpha_2$. Then
$\Phi^+=\{\pm\alpha_1+n\delta \mid n\in\dbZ\}\cup
\{\pm 2\alpha_1+(2n+1)\delta \mid n\in\dbZ\}$
(see \cite[Exercise~6.6]{Kac}).
For each $\alpha\in \Phi$ define an element $e_{\alpha}\in {\mathfrak sl}_3(k[t,t^{-1}])$
as follows:

\vskip .4cm
\renewcommand{\arraystretch}{1.5}
\setlength{\arrayrulewidth}{0.6pt}
\centerline{
\begin{tabular}{|c|c|}
\hline
$\alpha$ & $e_{\alpha}$\\
\hline
$\alpha_1+2n\delta$&$(e_{12}+e_{23})t^{2n}$\\
\hline
$\alpha_1+(2n+1)\delta$&$2(e_{12}-e_{23})t^{2n+1}$\\
\hline
$2\alpha_1+(2n+1)\delta$&$e_{13}t^{2n+1}$\\
\hline
$-\alpha_1+2n\delta$&$2(e_{21}+e_{32})t^{2n}$\\
\hline
$-\alpha_1+(2n+1)\delta$&$(e_{21}-e_{32})t^{2n+1}$\\
\hline
$-2\alpha_1+(2n+1)\delta$&$e_{31}t^{2n+1}$\\
\hline
\end{tabular}}
\vskip .4cm

\begin{proposition}\label{app}
Let $\Ghat=\Ghat(C)$, $\Uhat=\Uhat(C)$ and $\Gcal=\{g\in SL_3(k((t)))\mid J^{-1}g^* J=g^{-1}\}$. Then $\Ghat$
is isomorphic to $\Gcal$ (as a topological group) via the
map $\iota$ defined by
\[
	\iota: \chi_{\alpha}(u)\mapsto 1+(ue_{\alpha})+(ue_{\alpha})^2/2
	\ \text{ for } \
	\alpha\in\Phi
	\ \text{ and }\
	u\in k.
\]
Furthermore, $\iota(\Uhat)=\Ucal$ where
$\Ucal=\{g\in \Gcal\,\cap\, SL_3(k[[t]])\mid g \mbox{ is upper-unitriangular}\mod t\}$.
\end{proposition}

\begin{remark}
The expression $1+(ue_{\alpha})+(ue_{\alpha})^2/2$ should really be ``thought of'' as
$\exp(u e_{\alpha})$ since $e_{\alpha}^3=0$ for every $\alpha\in\Phi$.
\end{remark}

The proof of Proposition~\ref{app} will be given in Appendix A.
Henceforth we identify $\Ghat$ with $\Gcal$ and $\Uhat$ with $\Ucal$.
Before proceeding, we introduce some terminology. Let $M_3(k)$ denote
the space of $3\times 3$ matrices over $k$.

\begin{definition}
Let $g\in GL_3^1(k[[t]])$.  Write $g$ in the form $1+\sum_{i\geq 1} g_i t^i$
where $g_i\in M_3(k)$, and let $n$ be the smallest integer such that $g_n\neq 0$.
We will say that $g$ has \it{degree }\rm $n$ and write $\deg(g)=n$. The matrix $g_n$
will be called the \it{leading coefficient }\rm of $g$; we will write $\LC(g)=g_n$.
\end{definition}

Given a subgroup $H$ of $GL_3^1(k[[t]]])$ and $n\geq 1$, let
\[
    L_n(H)=\{\LC(g)\mid g\in H\mbox{ and
    }\deg(g)=n\}\cup\{0\} \, .
\]
Then it is easy to see that $L_n(H)$ is an $\mathbb F_p$-subspace of $M_3(k)$.
The following result is also straightforward.
\begin{lemma}
Let $S$ be a subgroup of $GL_3^1(k[[t]])$. For each $i\geq 1$ let
$S_i=S\cap GL_3^i(k[[t]]])$. Fix $n\in\dbN$, and let
$X\subset S$ be a set of elements of degree $n$. Then
$X$ generates $S_n$ modulo $S_{n+1}$ if and only if
the set $\{\LC(g)\mid g\in X\}$ spans $L_n(S)$.
\label{Liealg}
\QED
\end{lemma}

Now we return to the proof of Theorem~\ref{rank2}.
Let $P=\Uhat$. Define the filtration $\{P_n\}_{n=1}^{\infty}$ of $P$
as follows:
if $n\geq 2$, set $P_n=P\cap GL_3^n(k[[t]])$, and let $P_1$
be the set of matrices in $P\cap GL_3^1(k[[t]])$ whose
$(3,1)$-entry lies in $t^2 k[[t]]$.

Consider the following elements of $\Ghat$:
\[
	\{e_n^{(1)}(u), f_n^{(1)}(u), h_n(u)\mid n\in\Z, u\in k\}
	\ \text{ and } \
	\{e_n^{(2)}(u), f_n^{(2)}(u) \mid n\mbox{ is odd, } u\in k\},
\]
where
\begin{align*}
    e_n^{(1)}(u) &= \chi_{\alpha_1+n\delta}(u),
    &f_n^{(1)}(u) &= \chi_{-\alpha_1+n\delta}(u),
    \ \ \text{ and }\ \
    h_n(u)=[e_0^{(1)}(u), f_n^{(1)}(1)] \\
    e_n^{(2)}(u) &= \chi_{2\alpha_1+n\delta}(u),
    &f_n^{(2)}(u) &= \chi_{-2\alpha_1+n\delta}(u) .
\end{align*}
Let $E_n^{(i)}$, $F_n^{(i)}$ and $H_n$ be the subsets
$\{e_n^{(i)}(u)\}$, $\{f_n^{(i)}(u)\}$ and $\{h_n(u)\}$, respectively.

Now consider the subsets $\{Z_n\}_{n=1}^{\infty}$  of $P_1$ defined as follows:
\begin{eqnarray*}
    Z_1 &=& E_1^{(2)}\cup E_1^{(1)}\cup H_1 \cup F_1^{(1)},
    \qquad
    Z_{2n} = E_{2n}^{(1)}\cup F_{2n}^{(1)}\cup H_{2n}, \quad
    \text{ and } \\
    &&
    Z_{2n+1} = E_{2n+1}^{(1)}\cup F_{2n+1}^{(1)}\cup
    H_{2n+1}\cup E_{2n+1}^{(2)}\cup F_{2n+1}^{(2)}
\text{ for } n\geq 1.
\end{eqnarray*}

We claim that $Z_n$ generates $P_n$ modulo $P_{n+1}$ for each $n\geq 1$.
This result follows directly from Lemma~\ref{Liealg} applied with
$S=P_1$ and $\{S_n\}=\{P_n\}$.
Indeed, for each $n\geq 1$, define $L_n\subset M_3(k)$ as follows:
\[
L_n= \{g\in M_3(k)\mid g^t J = (-1)^n Jg\} \, 
\mbox{ where } g^t \mbox{ is the transposed of } g.
\]

It is clear from the definitions that
$L_n(P_1)\subseteq L_n(\Uhat)= L_n$.
On the other hand, direct computation shows that for $n\geq 2$,
all (non-identity) elements of $Z_n$ have degree $n$
and their leading coefficients span $L_n$, so
$L_n(P_1)=L_n$.
Similarly, one shows that the leading coefficients of elements
of $Z_1$ span $L_1(P_1)$.

In order to finish the proof by using Lemma~\ref{wellbeh}, we need suitable
commutation relations between the elements $\{e_n^{(i)}(u), f_n^{(i)}(u), h_n(u)\}$.
Once again, these are obtained by direct computation:
\begin{align*}
e_{2n}^{(1)}(u) &\equiv [e_{1}^{(2)}(1),f_{2n-1}^{(1)}(-u)] \mod P_{2n+1} \quad
& e_{2n+1}^{(1)}(u) &\equiv [e_{1}^{(2)}(1),f_{2n}^{(1)}(u)] \mod P_{2n+2} \\
f_{2n}^{(1)}(u) &\equiv [f_{1}^{(2)}(1),e_{2n-1}^{(1)}(u)] \mod P_{2n+1} \quad
&f_{2n+1}^{(1)}(u) &\equiv [f_{1}^{(2)}(1),e_{2n}^{(1)}(-u)] \mod P_{2n+2} \\
e_{2n+1}^{(2)}(u) &\equiv [e_{0}^{(1)}(1),e_{2n+1}^{(1)}(-u/4)] \mod
P_{2n+2} \quad
& f_{2n+1}^{(2)}(u) &\equiv [f_{1}^{(1)}(1),f_{2n}^{(1)}(-u/4)] \mod
P_{2n+2}.
\end{align*}
Finally, elements $\{h_n(u)\}_{n\geq 1}$ lie in $[P,P]$ by definition.

\appendix\section*{Appendix A: On explicit realization of twisted affine Kac-Moody groups}

\newtheorem*{propb}{Proposition~\ref{app}}

In this section we prove Proposition~\ref{app}. Recall that
\begin{equation} \label{CJ}
	C=
	\begin{pmatrix}
		2  & -4\\
		-1 &  2\\
	\end{pmatrix}
	\ \text{ and } \
	J= \begin{pmatrix}
		0&0&1\\
		0&-1&0\\
		1&0&0\\
	\end{pmatrix} .
\end{equation}

\begin{propb}
Let $\Ghat=\Ghat(C)$, $\Uhat=\Uhat(C)$ and $\Gcal=\{g\in SL_3(k((t)))\mid
J^{-1}g^* J=g^{-1}\}$. Then $\Ghat$
is isomorphic to $\Gcal$ (as a topological group) via the
map $\iota$ defined by
\begin{equation} \label{iota}
	\iota: \chi_{\alpha}(u)\mapsto 1+(ue_{\alpha})+(ue_{\alpha})^2/2
	\ \text{ for } \
	\alpha\in\Phi
	\ \text{ and } \
	u \in k.
\end{equation}
Furthermore, $\iota(\Uhat)=\Ucal$ where
$\Ucal=\{g\in \Gcal\,\cap\, SL_3(k[[t]])\mid g \text{ is
upper-unitriangular}\mod t\}$.
\end{propb}

\proof
We proceed in several steps.

\emph{Step 1: } First, we claim that (\ref{iota})
defines a (unique) homomorphism $\iota_0$ from the incomplete group $G=G(C)$ to $SL_3(k[t,t^{-1}])$.
This follows directly from the presentation of $G$ by generators and relators.

\emph{Step 2: } Let $K=\Ker\iota_0$ and let $Z$ be the kernel of the natural
map $G\to \Ghat$.
At this step we show that $K\subseteq Z$. Recall that $(\T_{build})$ denotes the building topology on $G$.
Let $(\T_{aux})$ be the topology on $\iota_0(G)$ given by the base $\{\iota(V)\}$
where $V$ runs over subgroups of $G$ open in $(\T_{build})$.
Let $\Gtil$ be the completion of $\iota_0(G)$ with respect to $(\T_{aux})$.
Clearly, there exists a continuous homomorphism $\eps:\Ghat\to\Gtil$. By \cite[Theorem 2.A.1]{Re2}, $\Ghat$
is topologically simple, whence $\eps$ is injective (since $\Gtil$ is clearly non-trivial).
Thus we conclude that $K\subseteq Z$. It follows immediately that the map $\iota_0: G\to\iota_0(G)$
canonically extends to an isomorphism of topological groups $\iota:\Ghat\to\Gtil$.

\emph{Step 3: } Consider the topology $(\T_{cong})$ on $\iota_0(G)$
induced from the congruence topology on $SL_3(k((t)))$, and let $\Gbar$
be the completion of $\iota_0(G)$ with respect to $(\T_{cong})$.
At this step we show that $\Gbar$ coincides with $\Gcal$.
It is clear that $\Gbar\subseteq \Gcal$ since $\iota_0(G)\subset \Gcal$ by construction
and $\Gcal$ is closed in the congruence topology.

Now we prove the reverse inclusion $\Gcal\subseteq \Gbar$.
The group $\Gcal$ is an isotropic simple algebraic group over the local field $k((t))$
and hence has a BN-pair $(\Bcal, \Ncal)$ given by Bruhat-Tits theory.
An explicit description of $(\Bcal, \Ncal)$ is given in \cite[1.15]{Ti}:
$\Bcal=\{g\in \Gcal\,\cap\, SL_3(k[[t]])\mid g \mbox{ is upper-triangular}\mod t\}$
and $\Ncal$ is the semi-direct product of
the group
\[
	D :=
	\left\{
	\begin{pmatrix}
		x&0&0\\
		0&1&0\\
		0&0&\sig(x)^{-1}\\
	\end{pmatrix}
	\ {\Bigg|} \
	x\in k((t))
	\right\}
\]
and the group of order $2$ generated by the matrix $J$ defined previously in
\eqref{CJ}. Note that $\Bcal=\iota(H) \Ucal $ where $H$ is the diagonal subgroup
of $G$.

Let $\Ubar$ be the closure of $\iota_0(U)$ in $\Gcal$.
Applying Lemma~\ref{wellbeh} with $P=\Ucal$ and \mbox{$\{P_i\}=\{P\,\cap\,
SL_3^i(k[[t]])\}$}
and arguing as in Section~\ref{sec:ranktwo},
we conclude that $\Ubar=\Ucal$.
From the explicit description of $\Ncal$, it is clear that $\Ncal=(\Ncal\cap
\Ucal) \iota_0(N)$.
Thus $\Gbar$ contains both $\Ucal$ and $\Ncal$.
Since $(\Bcal, \Ncal)$ is a BN-pair and $\Bcal\subset \iota_0(N) \Ucal$, it
follows that
$\Gcal$ is generated by $\Ucal$ and $\iota_0(N)$. Since $\Gbar$ contains
$\Ubar=\Ucal$ and $\iota_0(N)$,
we conclude that $\Gbar=\Gcal$.

\it{Step 4: }\rm Now we prove that the groups $\Gbar$ and $\Gtil$ are topologically isomorphic.
Equivalently, we will show that the topologies $(\T_{cong})$ and $(\T_{aux})$ on $\iota_0(G)$ coincide.
The inequality $(\T_{aux})\leq (\T_{cong})$ is clear.
This inequality also implies that there is a natural homomorphism $\eps_1:\Gbar\to\Gtil$.
Since $\Gbar=\Gcal$ is simple, $\eps_1$ is injective.

Now suppose that $(\T_{cong})$ is strictly stronger than $(\T_{aux})$.
Since both topologies $(\T_{cong})$ and $(\T_{aux})$ are countably
based, it follows that there is a subgroup $V\subset\iota_0(G)$, open in $(\T_{cong})$,
and a sequence $\{g_n\}$ in $\iota_0(G)$ such that

(a) $g_n$ converges to $1$ with respect to $(\T_{aux})$;

(b) $g_n\not\in V$ for all $n$.

Condition (a) implies that $g_n\in\iota_0(B)$ for all sufficiently large $n$.
Let $\Bbar$ be the closure of $\iota_0(B)$ in $\Gbar$. Clearly, $\Bbar$ is compact and countably based,
so there exists a subsequence $\{g_{n_k}\}$ which converges to some $g\in\Bbar$.
Since $\{g_{n_k}\}$ converges to $1$ with respect to $(\T_{aux})$, it follows
that $g$ lies in the kernel of $\eps_1:\Gbar\to\Gtil$. Since $\eps_1$
is injective, we conclude that $g=1$, contrary to condition (b).

\it{Step 5: }\rm
Combining steps 3 and 4, we get $\Gtil=\Gbar=\Gcal$, so the map $\iota:\Ghat\to\Gtil$
defined at the end of step 2 has the desired properties.
The equality $\iota(\Uhat)=\Ucal$ holds by construction.
\QED

\end{document}